\documentclass[11pt]{amsart}
\usepackage{amsmath,amscd,amssymb,amsfonts,amsthm,hyperref}

\usepackage{xcolor}  
\hypersetup{   colorlinks,   linkcolor={red!50!black},    citecolor={blue!70!black},   urlcolor={blue!80!black}}

\usepackage[cmtip,matrix,arrow]{xy}

\newtheorem{theorem}{Theorem}[section]

\newtheorem{proposition}[theorem]{Proposition}

\newtheorem{lemma}[theorem]{Lemma}
\theoremstyle{definition}    
\newtheorem{definition}[theorem]{Definition}
\theoremstyle{remark}

\newtheorem{remark}[theorem]{Remark}

\newtheorem{example}[theorem]{Example}
\newtheorem{examples}[theorem]{Examples}
\renewcommand{\AA}{\mathbb{A}}

\newcommand{\K}{\mathbb{K}}
\newcommand{\W}{\mathcal{W}}

\renewcommand{\L}{\mathcal{L}}
\renewcommand{\O}{\mathcal{O}}

\newcommand{\J}{\mathcal{J}}

\newcommand{\ca}{\mathcal}

\newcommand{\U}{\on{U}}
\newcommand{\E}{\ca{E}}

\newcommand{\N}{\mathbb{N}}
\newcommand{\R}{\mathbb{R}}
\newcommand{\C}{\mathbb{C}}

\newcommand{\Z}{\mathbb{Z}}

\newcommand\pt{\on{pt}}

\newcommand\lie[1]{\mathfrak{#1}}
\renewcommand{\k}{\lie{k}}
\newcommand{\h}{\lie{h}}
\newcommand{\g}{\lie{g}}

\renewcommand{\t}{\lie{t}}

\newcommand{\on}{\operatorname}

\newcommand{\Hom}{ \on{Hom}}

\newcommand\qu{/\kern-.7ex/} 

\newcommand{\lra}{\longrightarrow}

\renewcommand{\d}{{\mbox{d}}}
\newcommand{\ol}{\overline}

\newcommand{\f}{\frac}

\newcommand{\p}{\partial}

\newcommand{\ti}{\tilde}
\newcommand{\eeq}{\end{eqnarray*}}
\newcommand{\beq}{\begin{eqnarray*}}

\renewcommand{\H}{\ca{H}}
\renewcommand{\K}{\ca{K}}

\newcommand{\wh}{\widehat}
\newcommand{\lin}{{\on{lin}}}
\newcommand{\pol}{{\on{pol}}}

\newcommand{\mf}{\mathfrak}
\newcommand{\rra}{\rightrightarrows}
\renewcommand{\subset}{\subseteq}
\renewcommand{\supset}{\supseteq}
\newcommand{\I}{\ca{I}}

\newcommand{\wt}{\widetilde}

\newcommand{\su}{\mathsf{u}}
\newcommand{\sj}{\mathsf{j}}
\newcommand{\si}{\mathsf{i}}

\newcommand{\alg}{{\on{alg}}}

\begin{document}

\title{Differential Geometry of Weightings} 
\author{Yiannis Loizides}
\author{Eckhard Meinrenken}

\begin{abstract} We describe the notion of a \emph{weighting} along a submanifold $N\subset M$, 
and explore its differential-geometric implications. This includes a detailed discussion of weighted normal bundles, weighted deformation spaces, and weighted blow-ups. We give a description of weightings in terms of subbundles of higher tangent bundles, which leads us to notions of multiplicative weightings for Lie algebroids and Lie groupoids. 
\end{abstract}

\maketitle

\tableofcontents
\section{Introduction}
The technique of assigning weights to local coordinate functions is used in many areas of mathematics, such as singularity theory,  microlocal analysis, sub-Riemannian geometry, or algebraic geometry, under various terminologies. A  weighting along a submanifold $N\subset M$ may be described in  terms of a filtration 
\[ C^\infty_M=C^\infty_{M,(0)}\supset C^\infty_{M,(1)}\supset \cdots\]
of the sheaf of functions whose $i$-th term is  thought of as the sheaf of functions vanishing to \emph{weighted} 
order $i$ along $N$. See Section \ref{sec:weightings} for a precise definition. The idea for such `quasi-homogeneous structures' in terms of filtrations goes back to Melrose \cite{mel:cor}. 

As an example, letting $N=\{0\}$ be the origin in $M=\R^3$,  one can consider a weighting where the coordinate functions $x,y,z$ are assigned weights $1,2,3$ respectively. 
The ideal of functions of filtration degree four is then generated by the monomials $x^4,\,x^2 y,\,xz,\,yz,\,y^2,z^2$. 

A weighting along a submanifold $N\subset M$ determines a decreasing filtration of the normal bundle 
\[ \nu(M,N)=F_{-r}\supset \cdots \supset F_{-1}\supset F_0=0;\]
here $r$ is the \emph{order} of the weighting (an upper bound for the weights of local coordinate functions), 
and $F_{-i}$ is spanned by normal directions annihilating all functions of filtration degree $ i+1$. 
It turns out that weightings of order $r=2$ are completely determined by the subbundle $F=F_{-1}$. However, for $r>2$ knowledge of the  $F_{-i}$ does not suffice: intuitively, these subbundles  
only give first-order information whereas the weighting requires higher order information. This can be made precise using jet bundles. Letting $T_rM=J_0^r(\R,M)$ be the $r$-th tangent bundle, we will prove that an order $r$ weighting is equivalent to a certain subbundle
\[ Q\subset T_rM\]
along $N\subset M$. Theorem \ref{th:intrinsic} gives an intrinsic characterization of subbundles corresponding to weightings. 
The new description of weightings has some advantages; for example, one is led to notions of multiplicative weightings   
on Lie groupoids and infinitesimally multiplicative weightings on Lie algebroids.  
  
One of our main observations is that an order $r$ weighting along a submanifold determines a fiber bundle over $N$, which we call the \emph{weighted normal bundle}: 
\[ \nu_\W(M,N)\to N.\]
It admits an algebraic definition as the character spectrum of the associated graded algebra  of the algebra of smooth functions, or an alternative description as a subquotient $Q/\!\!\sim$ of the $r$-th tangent bundle $T_rM$. (For a trivial weighting, the latter is the usual definition  $\nu(M,N)=TM|_N/TN$ of the normal bundle.) 
The weighted normal bundle is not naturally a vector bundle, but it comes with an action of the multiplicative monoid $(\R,\cdot)$, making into a \emph{graded bundle} in the sense of Grabowski-Rotkiewicz \cite{gra:hig}.  It is 
non-canonically isomorphic to the associated graded bundle of $\nu(M,N)$, for the filtration described above.
In concrete examples, the weighted normal bundle often has additional structure as a bundle of nilpotent Lie groups, or as a bundle of homogeneous spaces of nilpotent Lie groups.

The properties of the ordinary normal bundle extend to the weighted setting. For example, a function $f\in C^\infty(M)$ of filtration degree $i$ canonically determines a function $f^{[i]}\in C^\infty(\nu_\W(M,N))$ on the weighted normal bundle,  homogeneous of degree $i$. Similarly, one has weighted homogeneous approximations of 
more general tensor fields, such as 
differential forms and vector fields. 
There is also a \emph{weighted deformation space} 
\[ \delta_\W(M,N)=\nu_\W(M,N)\sqcup (M\times \R^\times)\stackrel{\pi}{\lra} \R,\]
allowing to interpolate between a tensor field of given filtration degree with its weighted homogeneous approximation. 

As one important application of the deformation spaces, one obtains a simple coordinate-free construction of \emph{weighted blow-ups} of $M$ along $N$, extending many of the properties of ordinary blow-ups along submanifolds.  Recall that weighted blow-ups  in algebraic geometry have recently been used to great effect to give a short proof of Hironaka's desingularization theorem \cite{abr:fun,mcq:fas}. 

More generally, one can consider \emph{multi-weightings} for a finite collection of submanifolds $N_1,\ldots,N_d$. We shall require 
that these have \emph{clean intersection}, in the sense that near every $m\in M$ the intersection of submanifolds passing through $m$ is modeled by an intersection of coordinate subspaces. There are corresponding multi-weighted normal bundles $\nu_\W(M,N_1,\ldots,N_d)$ and deformation spaces 
\[ \delta_\W(M,N_1,\ldots,N_d) \stackrel{\pi}{\lra} \R^d.\]  For any multi-weighting, one obtains a total weighting along the intersection $N_1\cap\cdots\cap N_d$. If $d\ge 2$, the total weighting is  non-trivial even if the original multi-weighting is just given by the order of vanishing on the $N_\alpha$'s.   As a special case, every nested sequences of submanifolds $N_1\subset N_2\subset \cdots \subset N_r$ determines a weighting.  

Weightings, and the more general multi-weightings, appear in many contexts: We have already mentioned that a weighting of order $r=2$ amounts to the choice of a subbundle of the normal bundle.  As a simple example, given submanifolds $N_1,N_2\subset M$ with clean intersection $N$, the total weighting for the corresponding 
	trivial bi-weighting is given by  the subbundle $(TN_1|_N+TN_2|_N)/TN\subset \nu(M,N)$. 
	In symplectic geometry, given an isotropic submanifold $N$ of a symplectic manifold $(M,\omega)$, the symplectic normal bundle $TN^\omega/TN\subset \nu(M,N)$ defines a 2nd order weighting.  In \cite{me:eul}, this is used to prove a version of the Weinstein isotropic embedding theorem, where the `local model' does not involve choices.  
Collections of cleanly intersecting submanifolds, and the corresponding multi-weightings,  arise in the theory of manifolds with corners (e.g., \cite{mel:cor}), as well as in Lie theory  (e.g., \cite{beh:the}). Another source of examples are manifolds with \emph{Lie filtrations}, given by a filtration of the tangent bundle 
such that the induced filtration of sections is compatible with Lie brackets.  These so-called \emph{filtered manifolds} have been studied deeply in recent years, owing to their significance in the theory of hypo-elliptic operators. See 
\cite{choi:tan,hig:eul,moh:ind,moh:def,erp:ind,erp:tan,erp:gro}, for example. The  \emph{osculating tangent bundle}  
and \emph{osculating tangent groupoid} (using the terminology from \cite{erp:gro}) may be regarded as the 
weighted normal bundle and deformation space of the diagonal $M\subset M\times M$, for a suitable weighting 
defined by the Lie filtration. Similarly,  the normal bundles for filtered submanifolds, as 
in \cite{hig:eul}, are realized as weighted normal bundles. 

The full details of these applications to filtered manifolds, and a further generalization to singular Lie filtrations, will be given in a separate article. Forthcoming work of Daniel Hudson uses the weighted deformation spaces to develop a theory of 
weighted conormal distributions, complementing the work of Van Erp-Yuncken \cite{erp:tan} and Debord-Skandalis 
\cite{deb:lie}.  For future work, we envisage applications to normal forms for vector fields and 
Poisson structures in weighted settings (see e.g. \cite[Section 2.2]{duf:po} or \cite{lau:po}).

\bigskip\bigskip
\noindent{\bf Acknowledgements.} We thank Iakovos Androulidakis, Malte Behr, 
Nigel Higson, Daniel Hudson, Camille Laurent-Gengoux, Richard Melrose, Omar Mohsen, Brent Pym, Erik van Erp, and Robert Yuncken for discussions and helpful comments.
\bigskip\bigskip

\section{Weightings along submanifolds}\label{sec:weightings}
\subsection{Definitions}\label{subsec:def}
We start out with a definition of weightings on $n$-dimensional manifolds in terms of a local model. 


\begin{definition}
A \emph{weight sequence}  is a non-decreasing sequence of non-negative integers 
\begin{equation}\label{eq:weightsequence} 0\le w_1\le w_2\le \cdots \le w_n.\end{equation}
\end{definition}
The weight sequence  is  equivalently determined by the numbers $0\le k_0\le k_1\le \cdots\le n$, where 
\[ a\le k_i\ \Leftrightarrow\ w_a\le i.\]
That is, $k_i=\#\{a|\ w_a\le i\}$. An upper bound $r\in \N$ of the weight sequence we be called called its \emph{order}. 
(We don't insist that this is the least upper bound, hence, a weight sequence of order $r$ may also be regarded as a 
weight sequence of order $r'\ge r$.) By definition, $k_r=n$. 

Given an open subset  $U\subset \R^n$ and any $i\ge 0$, let 
$C^\infty(U)_{(i)}$ be the ideal generated by monomials 
\begin{equation}\label{eq:monomial} x^s=x_1^{s_1}\cdots x_n^{s_n},\ \ \ s=(s_1,\ldots,s_n)\end{equation}
with $s\cdot w\equiv\sum_{a=1}^n s_a w_a\ge i$. These ideals determine a 
filtration of the algebra of smooth functions
 \begin{equation} \label{eq:coordinatefiltration} C^\infty(U)=C^\infty(U)_{(0)}\supset C^\infty(U)_{(1)}\supset \cdots,\end{equation}
 compatible with the algebra structure. Some immediate observations:
	\begin{enumerate}
		\item In the definition of $C^\infty(U)_{(i)}$, we need only consider  monomials $x^s$ such that $s$ is 
		minimal relative to the condition $s\cdot w\ge i$; in particular, we may demand that  $s_a=0$ for $w_a=0$.
		\item If  $\R^{k_0}\cap U=\emptyset$, the filtration is trivial: $C^\infty(U)_{(i)}=C^\infty(U)$ for all $i$.
	\end{enumerate}

\begin{definition}
An \emph{order $r$ weighting}  on a $C^\infty$   manifold $M$ is a filtration of its	structure sheaf by ideals,
\begin{equation}\label{eq:sheaf} C^\infty_M=C^\infty_{M,(0)}\supset C^\infty_{M,(1)}\supset \cdots\end{equation}
with the property that every point of $M$ has an open neighborhood $U$, with coordinates $x_a$, such that 
the filtration of $C^\infty_M(U)=C^\infty(U)$ is given by \eqref{eq:coordinatefiltration}.   A \emph{morphism of weighted manifolds} is a smooth map $\phi\colon M\to M'$ such that the pullback map $\phi^*\colon C^\infty_{M'}\to C^\infty_M$ preserves filtrations. 
\end{definition}
We refer to the special coordinates in this definition as \emph{weighted coordinates}. Note that the coordinate function $x_a\in C^\infty(U)$ has filtration degree $w_a$. Observe also that an order $r$ weighting may be regarded as an order $r'$ weighting, for any $r'\ge r$. \smallskip
\begin{remark}
After completing a first version of this paper, we learned the concept of weightings had been introduced in lecture notes of Richard Melrose from 1996, under the name of \emph{quasi-homogeneous structure}.  See \cite[Proposition 1.15.1]{mel:cor}. The details of Melrose's approach were further developed in the recent Ph.D. thesis of Malte Behr \cite{beh:the}. 
\end{remark}


\begin{lemma}
	The ideal $C^\infty_{M,(1)}$ is the ideal sheaf $\I=\I_N$ of a closed submanifold $N\subset M$. A 	
	morphism of weighted manifolds $\phi\colon M\to M'$ takes $N$ into the corresponding closed submanifold $N'\subset M'$.  
\end{lemma}	
\begin{proof}
	Let $N\subset M$ be the closed subset of all points where the filtration \eqref{eq:sheaf} is non-trivial. 
	Given $m\in N$, choose 
	weighted coordinates $x_a$ on an open neighborhood $U$ of $m$. From the definition of the filtration in coordinates, 
	we see that $C^\infty(U)_{(1)}$ is the ideal of functions vanishing on $\R^{k_0}\cap U$. In particular, 
	the weighted coordinates give submanifold charts for $N$. The last claim is immediate from $\phi^* C^\infty_{M',(1)}\subset 
	C^\infty_{M,(1)}$. 
\end{proof}

If a closed submanifold $N\subset M$, with ideal sheaf $\I$,  is given in advance, we refer to a 
weighting with $C^\infty_{M,(1)}=\I$ as a \emph{weighting along $N$}, and to $(M,N)$ as a 
\emph{weighted manifold pair}.  

\begin{remark}\label{rem:closednogood}
The requirement that the submanifold $N$ be \emph{closed} is sometimes inconvenient. One can give a more flexible 
definition, for arbitrary  immersions $\si\colon N\to M$, as a filtration of the inverse image sheaf $\si^{-1}C^\infty_M$. 	Alternatively, given any $m\in N$, one can replace $N$ with a suitable neighborhood of $m$, and $M$ with a neighborhood of $\si(m)$, to reduce to the case that $N$ is closed. 
\end{remark}

From the local coordinate description of the filtration, we see that 
functions of filtration degree $i>r$ are sums of products of functions of filtration degrees 
$i_\nu\le r$, with $\sum i_\nu =i$. Hence, an order $r$ weighting of $M$ is  fully determined by the ideals $C^\infty_{M,(i)}$ with  $i\le r$:
\begin{equation}
i>r\ \Rightarrow \ C^\infty_{M,(i)}\subset \I^2.
\end{equation}
Recall that $\I/\I^2$ is the sheaf of sections of the conormal bundle
$\on{ann}(TN)\subset T^*M|_N$.

\begin{proposition}\label{prop:f}
	An order $r$ weighting of $(M,N)$ determines a filtration of the normal bundle $\nu(M,N)=TM|_N/TN$, 
	\begin{equation} \label{eq:normalbundlefiltration}
	\nu(M,N)=F_{-r}\supset \cdots F_{-1}\supset F_0=N,
	\end{equation}
	by subbundles $F_{-i}$ of dimension $k_i$, with the property that for all $i\ge 1$, the quotient
	\begin{equation}
		\label{eq:fquot}
	C^\infty_{M,(i)}/(C^\infty_{M,(i)}\cap \I^2)
	\end{equation}
	is the sheaf of sections of $\on{ann}(F_{-i+1})\subset \nu(M,N)^*$. 
\end{proposition}
\begin{proof}
The quotient  sheaf \eqref{eq:fquot} vanishes over $M-N$. At points $m\in N$, we may choose local weighted coordinates $x_a$ on an open neighborhood $U$ of $m$. The space of sections of \eqref{eq:fquot} over $N\cap U$ 
is then spanned by $\d x_a|_{N\cap U}$ 
 with  $w_a\ge i$, i.e., $a>k_{i-1}$.   
\end{proof}

The converse does not hold: A filtration of $C^\infty_M$, with $C^\infty_{M,(1)}=\I$ and 
$C^\infty_{M,(r+1)}\subset \I^2$, and  such that the quotients \eqref{eq:fquot} define subbundles of the conormal bundle, need not correspond to a weighting. 

\begin{example}
Let $A_i=C^\infty_{M,(i)}$ be given by 
	\[ A_1=\I,\ A_2=A_3=\I^2,\ A_4=A_5=\I^3,\ldots\]
	This is a multiplicative filtration with the property that the quotients  \eqref{eq:fquot} are zero for $i>1$, and equal to $\I/\I^2$ for	$i=1$. Clearly, this filtration does not correspond to a weighting. 
\end{example}	

The following proposition gives an additional condition ruling out such examples. Similar conditions are stated in \cite{beh:the,mel:cor}. 

\begin{proposition}
Let $N\subset M$ be a closed submanifold, with vanishing ideal $\I$, and let 
	\begin{equation}	\label{eq:sheaf} 
	C^\infty_M=C^\infty_{M,(0)}\supset C^\infty_{M,(1)}\supset \cdots
	\end{equation}
be a multiplicative filtration with $C^\infty_{M,(1)}=\I$ and $C^\infty_{M,(r+1)}\subset \I^2$. Then \eqref{eq:sheaf} 
defines a weighting of order $r$ along $N$ if and only if the following two conditions hold:
\begin{enumerate}
\item For $i\ge 1$, the quotient $C^\infty_{M,(i)}/(C^\infty_{M,(i)}\cap \I^2)$ is the sheaf of sections of a vector bundle.
\item For $i\ge 2$,
\[ C^\infty_{M,(i)}\cap \I^2=\sum_{1\le j<i} C^\infty_{M,(j)}\cdot C^\infty_{M,(i-j)}. \]
\end{enumerate}
\end{proposition}
\begin{proof}
Clearly, the filtration defined by a weighting has these properties. Conversely, if these conditions are satisfied, 
Condition (a) gives a filtration of the conormal bundle by subbundles
\[ \nu(M,N)^*=E_1\supset E_2 \supset \cdots \supset E_r \supset E_{r+1}=0, \]
where $E_i$ has $C^\infty_{M,(i)}/(C^\infty_{M,(i)}\cap \I^2)$ as its space of sections. 
Put $k_i=\dim M-\on{rank}(E_{i+1})$, so that $\dim N=k_0\le k_1\le \cdots \le k_r=n$, and define a weight sequence by letting $w_a=i$ for $k_{i-1}<a\le k_i$. 

Given a point $m\in N$, choose functions $x_a\in C^\infty(U)_{(r)}$ for $k_{r-1}<a\le n$, defined on an open neighborhood $U$ of $m$, such that their images under the quotient map  define a frame for $E_r$ over $U\cap N$. Taking $U$ smaller if needed, extend to a collection of functions $x_a\in C^\infty(U)_{(r-1)}$ for $k_{r-2}<a\le n$  such that their images under the quotient map  define a frame for $E_{r-1}$ over $U\cap N$. Proceeding in this way, we obtain local coordinates $x_1,\ldots,x_n$ near $m$,  with the property  
$x_a\in  C^\infty(U)_{(i)}$ for $k_{i-1}<a\le k_i$. Define a weighting on $U$ along $N\cap U$ where the coordinate function 
$x_a$ has weight $w_a$.  

Let $C^\infty_U=B_0 \supset B_1 \supset B_2\supset \cdots$ be the filtration defined by this weighting, and let 
$A_i=C^\infty_{U,(i)}$. We shall argue by induction that 
$A_i=B_i$ for all $i$. 
By assumption, $A_1=B_1=\I|_U$.  Suppose $i>1$ and $B_j=A_j$ for $j<i$. By condition (b) and the inductive hypothesis,
\[ B_i \cap \I^2=\sum_{1\le j< i} B_j\cdot B_{i-j}=\sum_{1\le j<i} A_j\cdot A_{i-j}=A_i \cap \I^2.\]
Note that $B_i$ is generated by this intersection, together with the coordinate functions $x_a$ with $k_{i-1}<a\le k_i$. Hence $B_i\subset A_i$. To prove equality, it suffices to show that the induced map
$B_i/(B_i\cap \I^2)\to A_i/(A_i\cap \I^2)$ is an isomorphism. But this is clear since  
 both are identified with the sheaf of sections of $E_i|_U$.  
\end{proof}

\subsection{Examples, basic constructions}
\begin{examples}\label{ex:ex}
	 Let $N\subset M$ be a closed submanifold. 
	\begin{enumerate}
		\item \label{it:0} A weighting of order $r=1$ along $N$ is the \emph{trivial weighting} along $N$, given by 
		\[ C^\infty_{M,(i)}=\I^i,\] the ideal of functions vanishing to order $i$ on $N$ (in the usual sense).
		\item \label{it:a1}  A weighting of order $r=2$ along $N$ is fully determined by the subbundle  $F=F_{-1}\subset \nu(M,N)$: 
	 Let  $\J\subset \I$ be the functions $f$ such that $f|_N=0$ and $\d f|_{\ti{F}}=0$, where 
	 $\ti{F}\subset TM|_N$ is the pre-image  of $F$. Then \cite{me:eul}
		\[ C^\infty_{M,(2i+1)}=\I\ca{J}^i,\ \ C^\infty_{M,(2i)}=\ca{J}^i.\]
		\item \label{it:nested}	
		Consider a nested sequence of submanifolds of $M$ 
		\[ M=N_{r+1}\supset N_r\supset \cdots \supset N_1=N,\]
		and let $\I_i\subset \I$  be the vanishing ideal of $N_i$. Then 
		\[ C^\infty_{M,(i)}=\sum_{\ell\ge 1}\sum_{i_1+\ldots+i_\ell= i} \I_{i_1} \cdots \I_{i_\ell}\] 
		for $i\ge 1$ defines an order $r$ weighting  along $N$. Near any point $m\in N$ we may choose  
		coordinates $x_1,\ldots,x_n$, that are submanifold coordinates for each of the $N_i$'s; these coordinates  will serve as weighted coordinates.
		The corresponding filtration \eqref{eq:normalbundlefiltration} 
			of the normal bundle is given as $ F_{-i}=TN_{i+1}|_N/TN$. 
\item\label{it:clean}
 Let $N_1,N_2\subset M$ be closed submanifolds with clean intersection. That is, $N=N_1\cap N_2$ is a submanifold, with 
 $TN=TN_1\cap TN_2$.  Then $C^\infty_M$ acquires a bi-filtration by ideals $\I_{(i_1,i_2)}$, whose  local sections are the functions vanishing to order $i_1$ on $N_1$ and to order $i_2$ on $N_2$. Let 
\[ C^\infty_{M,(i)}=\sum_{i_1+i_2=i} \I_{(i_1,i_2)}\]
Using local coordinates adapted to the clean intersection (see \cite{pik:the}), one verifies that this is a weighting of order $r=2$ along $N=N_1\cap N_2$. The corresponding subbundle $F\subset \nu(M,N)$ is $(TN_1|_N+TN_2|_N)/TN$; this subbundle, and hence the weighting, is trivial if and only if the intersection is transverse. This example generalizes to finite collections of closed submanifolds  $N_1,\ldots,N_r\subset M$ with clean intersection (see Section \ref{subsec:multiweightings} below).  
\item (\emph{Linear weightings}.) \label{it:linweighting}
Let $V\to M$ be a vector bundle, and $W\to N$ a subbundle along $N\subset M$.  
We adopt the Grabowski-Rotkiewicz approach \cite{gra:hig} to vector bundles in terms of their scalar multiplication; thus, $W$ is characterized as a submanifold invariant under scalar multiplication. A weighting of 
$V$ along $W$  will be called \emph{linear} if 
the corresponding filtration of $C^\infty_V$ is $(\R,\cdot)$-invariant. It induces a unique weighting of $M$ along $N$, in such a way that the bundle projection and inclusion of units are morphisms of weighted manifolds. 
\item (\emph{Lie filtrations}.) A \emph{Lie filtration} of a manifold $M$ is a filtration of the tangent bundle by subbundles, 
\[ TM=H_{(-r)}\supset \cdots H_{(-1)}\supset H_{(0)}=0_M\]
such that the resulting filtration of vector fields preserves brackets. Such a Lie filtration canonically determines a weighting
of $M\times M$ along the diagonal $M$. Given a `filtered submanifold' $N\subset M$ \cite[Section 7]{hig:eul}, one 
obtains a weighting of $M$ along $N$.  These examples will be discussed in detail in \cite{loi:lie}, with further generalizations to singular Lie filtrations. 
	\end{enumerate}
\end{examples}

Some basic constructions with weightings:
\begin{enumerate}
	\item {\bf Products.} Given weightings of $(M,N)$ and $(M',N')$, 
	we obtain a weighting of the product  $(M\times M',N\times N')$, with $C^\infty_{M\times M',(i)}$ the ideal generated by the sum of all $C^\infty_{M,(i_1)}\otimes C^\infty_{M',(i_2)}$ with $i_1+i_2=i$. 
	\item{\bf Pullbacks.} Suppose we are given a weighting along $N\subset M$, and let $\varphi\colon M'\to M$ be a smooth map transverse to $N$. 
	Then there is a canonically induced weighting of $M'$ along $N'=\varphi^{-1}(N)$, 
	with $C^\infty_{M',(i)}$ the ideal generated by $\varphi^*C^\infty_{M,(i)}$. (To see that this is a weighting, let $m'\in N'$, and choose weighted coordinates $x_a$ near $m=\varphi(m')$. The pullbacks $\varphi^*x_a$ of coordinate functions with $a>k_0$ may be completed to a coordinate system near $m'$, and these will then serve as weighted coordinates near $m'$.) As a special case, given a submanifold $\Sigma\subset M$ transverse to $N$, 
	one obtains a weighting of $\Sigma$ along $\Sigma\cap N$. 
\end{enumerate}

\subsection{Filtrations on forms and vector fields}\label{subsec:filtrationonforms}
The filtration of $C^\infty_M$ extends uniquely to a filtration of the sheaf $\Omega^\bullet_M$ of differential forms, 
\[ \Omega^\bullet_M=\Omega^\bullet_{M,(0)}\supset \Omega^\bullet_{M,(1)}\supset \cdots\]
compatible with the algebra structure, in such a way that the de Rham differential preserves the filtration degree. 
In weighted coordinates on $U\subset M$,  the space $\Omega^q(U)_{(i)}$ consists of forms 
\[ \sum_{a_1<\cdots<a_q} f_{a_1\cdots a_q}\d x_{a_1}\wedge \cdots \d x_{a_q}\]
such that  $f_{a_1\cdots a_q}$ has filtration degree 
$i-w_{a_1}-\cdots-w_{a_q}$. 

Similarly, the sheaf of vector fields $\mf{X}_M$ inherits a filtration by $C^\infty_M$-submodules, 
where a vector field $X$ has filtration degree $i$ if the Lie derivative $\L_X$ on functions raises the filtration degree by $i$. This filtration starts in degree $-r$:
\[ \mf{X}_M=\mf{X}_{M,(-r)}\supset \mf{X}_{M,(-r+1)}\supset \cdots\]
In weighted coordinates  on $U\subset M$, the space 
$\mf{X}(U)_{(i)}$ consists of vector fields 
\[ \sum_a f_a\f{\p}{\p x_a}\]
such that $f_a$ has filtration degree $i+w_a$.  
These filtrations on vector fields and forms are compatible with all the operations from Cartan's calculus
(brackets, contractions, Lie derivatives, and so on).  

The sheaf $\mf{X}_{M,(0)}$ is closed under brackets. Its local sections are the vector fields 
preserving the filtration of $C^\infty_M$, and so are the  \emph{infinitesimal automorphisms} of the weighting.  
Using weighted coordinates on $U\subset M$, we see that vector fields in $\mf{X}(U)_{(0)}$ are tangent to $N\cap U$, and that the restriction map $\mf{X}(U)_{(0)}\to \mf{X}(N\cap U)$ is surjective. Hence, we obtain a surjective sheaf map $\mf{X}_{M,(0)}\to \mf{X}_N$.

\subsection{Other settings}
The definition of weightings was formulated in such a way that it extends to various other kinds of structure sheaves. 
In particular, for a complex manifold $M$ one defines a \emph{holomorphic weighting} in terms of a decreasing filtration of the sheaf $\O_M$ of holomorphic functions, given by the local model from Section \ref{subsec:def}
in  holomorphic coordinates. The degree $1$ part of such a filtration is the vanishing ideal of a complex submanifold $N$. The filtration of the structure sheaf determines filtrations on the sheaves of holomorphic differential forms and holomorphic vector fields.
\medskip

For much of the present paper, we will restrict the discussion to the $C^\infty$ category.  Having partitions of unity at our disposal, we may avoid the use of sheaves, and simply work with the filtration of global functions, 
\[	C^\infty(M)=C^\infty(M)_{(0)}\supset C^\infty(M)_{(1)}\supset \cdots.\]
In Section \ref{subsec:holo}, we will make a few more comments on the extension of our results to the holomorphic context.

\section{Graded bundles (after Grabowski-Rotkiewicz) }\label{sec:weightingsgraded}
In Section \ref{sec:weightednormalbundle}, we will define the weighted normal bundle for a given weighting along a submanifold. We shall see that the weighted normal bundle is not naturally a vector bundle; rather, it is a graded bundle in the sense of Grabowski-Rotkiewicz \cite{gra:gra}. To prepare, let us review the definition and basic properties of graded bundles, following \cite{bru:lin,bru:int,gra:dua2,gra:gra}.

\subsection{Definition and examples}
A graded bundle may be defined as a special case of a non-negatively graded supermanifold \cite{cat:sup,meh:sup} in which no fermionic variables are present. 
In this paper, we will work with a simple definition of graded bundles due to Grabowski-Rotkiewicz \cite{gra:gra}, solely in terms of scalar multiplications. 

\begin{definition}\cite{gra:gra} A \emph{graded bundle}  is a manifold $E$ with a smooth action 
	\begin{equation}\label{eq:kappaaction} \kappa\colon \R\times E\to E,\ (t,x)\mapsto \kappa_t(x)\end{equation}
	of the monoid $(\R,\cdot)$ of real numbers, $\kappa_{t_1t_2}=\kappa_{t_1}\circ \kappa_{t_2},\ \kappa_1=\on{id}_E$.	The  \emph{Euler vector field} of a graded bundle $E$ is the vector field $\E$ with flow $s\mapsto \kappa_{\exp(-s)}$.  A \emph{morphism of graded bundles} $\varphi\colon E'\to E$  is a smooth map  intertwining the $(\R,\cdot)$-actions. 
	
\end{definition}

For any graded bundle, the map $\kappa_0\colon E\to E$ is a smooth projection, hence its range  is a closed submanifold $N\subset E$. In fact, the projection $\kappa_0$ makes $E$ into a locally trivial fiber bundle over $N$.  A \emph{graded subbundle} is an $(\R,\cdot)$-invariant submanifold $E'\subset E$; it is a graded bundle in its own right, with the inclusion map a morphism of graded bundles. 	

\begin{examples}\label{ex:gradedbundles}
	\begin{enumerate}
		\item\label{it:a} A vector bundle is a graded bundle with the special property that the map $E\to TE,\ x\mapsto \f{d}{d t}|_{t=0}\kappa_t(x)$ is injective. Conversely, if this condition is satisfied, then $\kappa_t$ is the scalar multiplication for a unique vector bundle structure on $E$. 	 \cite{gra:gra}. 
		\item\label{it:b} 	Every negatively graded vector bundle 
		\[ W=W^{-1}\oplus \cdots \oplus W^{-r}\to N\]
		is a graded bundle, by letting $\kappa_t$
		be  multiplication by $t^i$ on $W^{-i}$. The corresponding Euler vector field is $\E=\sum_{i=1}^r i \E_i$, where $\E_i$ is the usual Euler vector field of the vector bundle $W^{-i}$.   	
		\item\label{it:c} Let $\g=\bigoplus_{i=1}^r \g^{-i}$ be a negatively graded Lie algebra.
		Then $\g$ is nilpotent, and the exponential map gives a global diffeomorphism with the 
		corresponding simply connected nilpotent Lie group $G$. 
		The $(\R,\cdot)$-action on $\g$ (as in \eqref{it:b}) is by Lie algebra morphisms, hence it exponentiates
		to an action by Lie group morphisms, making $G\to \pt$ into a \emph{graded Lie group}. These have appeared in the literature under the name of \emph{homogeneous Lie group} (see the monograph \cite{fis:qua}); if  
		$\g^{-1}$ generates $\g$ as a Lie algebra, they are also called \emph{Carnot groups}. 
		Given a graded Lie subalgebra 
		$\h\subset \g$, with corresponding subgroup $H\subset G$, the homogeneous space $G/H\to \pt$ becomes a graded bundle. 
		\item\label{it:d} 
		More generally, one can consider  \emph{graded Lie groupoids}: that is, Lie groupoids with an action of $(\R,\cdot)$ by Lie groupoid morphisms. (For a vector bundle, regarded as a Lie groupoid, this recovers the notion of graded vector bundle.) 
		\item \label{it:e}	The \emph{bundle of $r$-velocities} or 
		\emph{$r$-th tangent bundle}  of a manifold $M$ was introduced by Ehresmann \cite{ehr:pro} as the bundle of $r$-jets of curves, 
		\[ T_rM=J_0^r(\R,M)=\{\sj^r_0(\gamma)|\ \gamma\colon \R\to M\}.\]
		It is a graded bundle $T_rM\to M$, with the $(\R,\cdot)$-action coming from the action on the domain of 
		such paths: $\kappa_t(\sj^r_0(\gamma))=\sj^r_0(\gamma_t)$ with $\gamma_t(s)=\gamma(ts)$. 
		\item\label{it:f} According to Grabowski-Rotkiewicz \cite{gra:hig}, \emph{double vector bundles} may be defined as manifolds $D$ with two vector bundle structures $D\to A$ (`horizontal') 
		and $D\to B$ (`vertical'), in such a way that the corresponding scalar multiplications $\kappa_t^h,\kappa_t^v$ 
		commute. In this case, $A,B$ are themselves vector bundles, over a common base $N$. Double vector bundles become 
		graded bundles of order $2$ for the scalar multiplication $\kappa_t=\kappa_t^h\circ \kappa_t^v$.  More generally, one can consider multi-vector bundles. 
		
	\end{enumerate}	
\end{examples}

For every graded bundle $E\to N$, there is a canonical weighting 
of order $r$ along $N$ given by 
\begin{equation}\label{eq:weightingsgradedbundles} 	C^\infty(E)_{(i)}=\{f\in C^\infty(E)|\ \kappa_t^*f=O(t^i)\}.\end{equation}
The graded bundle coordinates on $E$ serve as weighted coordinates.


\subsection{Linear approximation}
Given a graded bundle $E\to N$, the normal bundle $\nu(E,N)$ becomes a graded bundle, with scalar multiplication $\nu(\kappa_t)$ obtained by applying the normal bundle functor to $\kappa_t\colon E\to E$. This 
\emph{linear approximation}  $E_\lin=\nu(E,N)$
is naturally a graded  \emph{vector} bundle, 
\[ E_\lin=\bigoplus_{i\ge 1} E_\lin^{-i},\]
with  $E_\lin^{-i}$ the sub-vector bundle on which 
$\nu(\kappa_t)$ acts as scalar multiplication by $t^i$. 
(Given a linear $(\R,\cdot)$-action on a real vector space $\mathsf{E}$, the complexified $(\C,\cdot)$-action on $\mathsf{E}\otimes_\R \C$ restricts to a $\U(1)$-action; the fact that it 
extends to a $(\C,\cdot)$-action means that negative weights cannot occur.) Morphisms of graded bundles 
$E\to E'$ induce morphisms of graded vector bundles $E_\lin\to E'_\lin$. 
We say that the graded bundle $E$ has \emph{order} $r$ if $E^{-i}_\lin=0$ for $i>r$. 

\begin{examples}\label{ex:linapp} The linear approximations for some of our examples are as follows:
	\begin{enumerate}
		\item For a graded vector bundle $W$, we have a canonical identification $W_\lin=W$.
		\item For a graded Lie group  $G$,  we recover the Lie algebra $G_\lin=\g$. The linear approximation of group multiplication $G\times G\to G$ is the addition $\g\times \g\to \g$. 
		\item For the $r$-th tangent bundle, we have that $(T_rM)_{\on{lin}}={TM\oplus \cdots \oplus TM}$. See Section \ref{subsubsec:action} below.
		\item \label{it:lind}
		For a double vector bundle $D$ with side bundles $A,\,B$ over a base $N$, let $\on{core}(D)\subset D$ be the subset on which $\kappa_t^h,\kappa_t^v$ coincide. This is a vector bundle over $N$, and one finds that  
		\[ D_\lin=(A\oplus B)\oplus \on{core}(D).\]
		with $A\oplus B$ in degree $-1$ and $\on{core}(D)$ in degree $-2$. 
	\end{enumerate}
\end{examples}

Grabowski-Rotkiewicz  established the following result, which may be seen as an analogue of 
Batchelor's theorem \cite{bat:str}  for supermanifolds:
\begin{theorem} \cite{gra:gra} \label{th:gr} For every graded bundle $E\to N$ there exists a 
	(non-canonical) 
	isomorphism of graded bundles
	\begin{equation} \label{eq:batch} E_{\on{lin}}\cong E\end{equation}
	whose linear approximation is the identity. 
\end{theorem}
We refer to any such isomorphism as a \emph{linearization} of $E$. When proving facts about graded bundles, it is often helpful to choose such a linearization to reduce to the case of a graded vector bundle.

\begin{remark}	
	There is an obvious version of graded bundles in the holomorphic or analytic categories. (See \cite{joz:monoid} for some details.) 
	In these settings, linearizations exists near any given $m\in N$, but need not exist globally. 	
\end{remark}


\subsection{Graded bundle coordinates}
Let $\pi\colon E\to N$ be a graded bundle of order $r$, with $\dim E=n$. Put $k_0=\dim N$, $k_i=\dim E_\lin^{-i}$ for $i>0$, and let $w_1\le \cdots \le w_n=r$   be the corresponding weight sequence,  given by $w_a=i$ for $k_{i}\le a<k_{i+1}$. 
Coordinates 
\[ x_a,\ \ a=1,\ldots,\dim E\]
defined on $\pi^{-1}(O)$ for open subsets  $O\subset N$, are called \emph{graded bundle coordinates} if 
$x_a$ for $k_i\le a< k_{i+1}$ are homogeneous of degree $i$. To construct such coordinates, choose a linearization, and pick  vector bundle coordinates for $E_\lin$ compatible with the grading. In graded bundle coordinates, the Euler vector field of $E\to N$ is given by 
\[ \E=\sum_a w_a x_a\f{\p}{\p x_a}. \]


\subsection{Polynomial functions, forms and vector fields}\label{subsec:polfunctions}
Given a graded bundle $E\to N$,  with Euler vector field $\E$, let $C^\infty(E)^{[k]}$ be the space of 
smooth functions $f$ that are \emph{homogeneous of degree $k$}, i.e. $\L_\E f=k\ f$. These spaces are non-zero only if $k$ is a nonnegative integer, and their direct sum
is a graded algebra $C^\infty_\pol(E)$, called the \emph{polynomial functions} on $E$. As shown in  \cite{gra:dua2}, 
one recovers $E$ from this algebra as the character spectrum 
\begin{equation}\label{eq:recovery} E=\on{Hom}_{\on{alg}}(C^\infty_\pol(E),\R);\end{equation}
here the $(\R,\cdot)$-action on algebra morphisms is given by duality with the $(\R,\cdot)$-action on polynomials. 
Similarly, we can consider vector fields or differential forms of homogeneity $k$, making 
\[ \mf{X}_\pol(E)=\bigoplus_{k\ge -r} \mf{X}(E)^{[k]},\ \ \ \ 
\Omega_\pol(E)=\bigoplus_{k\ge 0} \Omega(E)^{[k]}
\]
into a graded Lie algebra and a bigraded algebra, respectively. The gradings are compatible with the 
$C^\infty_\pol(E)$-module structure, as well as the operations from Cartan's calculus
(Lie derivatives, Lie brackets, contractions, differentials). Note that the grading of $\mf{X}_\pol(E)$ starts in degree $-r$; 
the part of strictly negative degree is a nilpotent Lie subalgebra 
\begin{equation}\label{eq:negative} \mf{X}_\pol(E)^{-}=\bigoplus_{1\le i\le r} \mf{X}(E)^{[-i]}.\end{equation}
The vector fields in this subalgebra annihilate $C^\infty(E)^{[0]}\cong C^\infty(M)$, hence 
they are all vertical. We have $\mf{X}(E)^{[-r]}\cong \Gamma(E_\lin^{-r})$ via restriction of vector fields to $M\subset E$, and identifying $TE|_M=TM\oplus E_\lin$. Analyzing fiberwise (and using graded bundle coordinates), one finds that 
the vector fields in \eqref{eq:negative} are all complete. Hence, this subalgebra exponentiates to a nilpotent group of fiber-preserving diffeomorphisms of $E$.  Observe also that 
\[ \mf{X}(E)^{[0]}=\mf{aut}_{\on{gr}}(E)\]
is the  Lie algebra of infinitesimal graded bundle automorphisms of $E$, given by $\kappa_t$-invariant vector fields. 
%

\section{Weighted normal bundle}\label{sec:weightednormalbundle}
The normal bundle of a closed submanifold can be defined in several equivalent ways: (i) as the character spectrum of the associated 
graded algebra of $C^\infty(M)$, using the filtration  by the powers of the vanishing ideal $\I$ of $N$, (ii) as a subquotient of the tangent bundle
$\nu(M,N)=TM|_N/TN$, (iii) in terms of its space of sections, $\Gamma(\nu(M,N))=\mf{X}(M)/\mf{X}(M,N)$ where 
$\mf{X}(M,N)$ is the subspace of vector fields tangent to $N$. Each of these descriptions admits a generalization to the weighted setting. In this section we follow approach (i), and also give a coordinate description. Later, we shall explain versions of (ii) and (iii). 
\subsection{Definition of the weighted normal bundle}
Let $(M,N)$ be a manifold pair with a weighting of order $r$. Consider the corresponding filtration of the algebra of smooth functions, and let $\on{gr}(C^\infty(M))$ be the associated graded algebra, with graded components 
\begin{equation}\label{eq:gradedcomponents}
C^\infty(M)_{(i)}/C^\infty(M)_{(i+1)}.\end{equation}
%
For $f\in C^\infty(M)_{(i)}$, its image in \eqref{eq:gradedcomponents} will be denoted $f^{[i]}$, and is called the ($i$-th) \emph{homogeneous approximation}  of $f$. 
If $i=0$, this is the restriction $f|_N\in C^\infty(N)$.
\begin{definition}[Weighted normal bundle] 
	\label{def:algebraic}
	The \emph{weighted normal bundle} for an order $r$ weighting along $N\subset M$ is the 
	character spectrum, 
	\[ \nu_\W(M,N)=\Hom_\alg(\on{gr}(C^\infty(M)),\R).\]
\end{definition}
By definition, $f^{[i]}$ for $f\in C^\infty(M)_{(i)}$ is a function on the weighted normal bundle. The 
smooth structure  is uniquely specified by the requirement that all such functions are smooth. More precisely: 

\begin{theorem}
	\begin{enumerate}
		\item The weighted normal bundle $\nu_\W(M,N)$ has a unique structure as a graded bundle over $N$, of dimension equal to that of $M$, in such a way that
		\[ C^\infty_\pol(\nu_\W(M,N))=\on{gr}(C^\infty(M)).\]
	   as graded algebras.
		\item  Given weighted coordinates $x_1,\ldots,x_n$ on $U\subset M$, the homogeneous approximations $(x_a)^{[w_a]}$ for $a=1,\ldots,n$ serve as graded bundle coordinates on $\nu_\W(M,N)|_{U\cap N}$. 		
	\end{enumerate}
\end{theorem}
\begin{proof}
This may be proved be using the general criterion in Haj-Higson \cite{hig:eul}, or by direct construction of coordinates, as follows. 
	Let $p\colon \nu_\W(M,N)\to N$ be the projection induced on the algebra side by the inclusion $C^\infty(N)\to\on{gr}(C^\infty(M))$  as the degree $0$ summand. For any open subset $U\subset M$, we have that 
	$p^{-1}(U)=\nu_\W(U,U\cap N)$. 
	Given weighted coordinates $x_1,\ldots,x_n$ on $U\subset M$, let  $y_a=(x_a)^{[w_a]}$ be their homogeneous approximations. 
	For $a\le k_0$ we have $y_a=x_a|_{U\cap N}$, hence $y_1,\ldots,y_{k_0}$ serve as coordinates on $U_1=U\cap N$. Let 
	$\ol{U}_1\subset \R^{k_0}$ be the image of this coordinate map.   The graded algebra 
	 $\on{gr}(C^\infty(U))$ is generated, as an algebra over its degree $0$ summand $C^\infty(U\cap N)$, by its elements
	$y_a$ for $a>k_0$. Hence, 
\[ (y_1,\ldots,y_n)\colon p^{-1}(U)\to \R^n\]
is a bijection onto an open subset of the form $\ol{U}_1\times \R^{n-k_0}$. 	We will take these as graded bundle charts. 
In particular, this specifies a $C^\infty$-structure on $p^{-1}(U)$. Let us show that if $f\in C^\infty(U)_{(i)}$, then $f^{[i]}\colon p^{-1}(U)\to \R$ is smooth. It suffices to check for functions of the form  $f(x_1,\ldots,x_n)=\chi(x_1,\ldots,x_n) x_1^{s_1}\cdots x_n^{s_n}$ with 
$\chi\in C^\infty(U)$ and $s\cdot w=\sum s_a w_a=i$. For such a function, 
\[ f^{[i]}(y_1,\ldots,y_n)=\chi(y_1,\ldots,y_{k_0},0,\ldots,0)\,y_1^{s_1}\cdots y_n^{s_n},\]
which is evidently smooth, and is homogeneous of degree $i$.   
This also show that for a different set $x_1',\ldots,x_n'\colon U'\to \R$ of weighted coordinates, the 
`transition function' is smooth: taking $x_a'$ as $f\in C^\infty(U\cap U')_{(w_a)}$, we see that $y_a'$ is smooth as a function of 
$y_1,\ldots,y_n$. This completes the construction of a $C^\infty$-structure for which  the functions $f^{[i]}$ with $f\in C^\infty(M)_{(i)}$ are all smooth. The uniqueness assertion is clear. 
\end{proof}

Note that for a weighted coordinate chart, the algebra $\on{gr}(C^\infty(U))$ coincides with the corresponding algebra 
for the trivial weighting on $U\cap N$ -- only the grading is different. This shows that locally, the weighted normal bundle 
is isomorphic to the usual normal bundle. However, this identification depends on the choice of coordinates, and the 
transition functions for the two bundles are different.

\begin{example}	
		Suppose $x,y,z$ are weighted coordinates  of weights $0,1,3$ on $\R^3$, and consider the coordinate change (defined near the origin) 
	\[ x'=\sin(x)\exp(yz),\ y'=y\exp(xy),\ z'=3z+\sin^3(xy).\]
	The corresponding coordinate change on the weighted normal bundle reads as
	\[ (x')^{[0]}=\sin(x^{[0]}),\ (y')^{[1]}=y^{[1]},\ \ (z')^{[3]}=3z^{[3]}+(x^{[0]})^3(y^{[1]})^3.\]
	The coordinate change for the usual normal bundle is given by a similar formula, omitting the `non-linear' term $(x^{[0]})^3(y^{[1]})^3$. 
\end{example}
The construction of the weighted normal bundle is functorial: 	Any morphism of weighted manifold pairs $\varphi\colon (M,N)\to (M',N')$ 
gives a morphism of	filtered algebras $C^\infty(M')\to C^\infty(M)$, hence of the associated graded algebras. 
It therefore determines  a map of  weighted normal bundles, 
$\nu(\varphi)\colon \nu_\W(M,N)\to \nu_{\W'}(M',N')$. To check that this map is smooth, it suffices to show that the pullback of smooth functions on $\nu_{\W}(M',N')$ is again smooth. Let 
$f\in C^\infty(M')_{(i)}$, with corresponding homogeneous function $f^{[i]}$ on $\nu_{\W}(M',N')$. 
Then $(\nu(\varphi))^*f^{[i]}=(\varphi^*f)^{[i]}$ is smooth, as required.

\subsection{Properties of the weighted normal bundle}
Recall that the weighting along $N\subset M$ determines a filtration \eqref{eq:normalbundlefiltration} of the (usual) normal bundle $\nu(M,N)$ by subbundles $F_{-i}$.  The associated graded bundle has its components in negative degrees, 
$\on{gr}(\nu(M,N))=\bigoplus_i F_{-i}/F_{-i+1}$. 
%
%
\begin{proposition}
	There is a canonical identification of graded vector bundles over $N$, 
	\[ 	\nu_\W(M,N)_\lin=\on{gr}(\nu(M,N)).\]
Here $\nu_\W(M,N)_\lin$ is the linear approximation of the graded bundle $\nu_\W(M,N)$. 
\end{proposition} 
\begin{proof}
	The space of polynomial functions of degree $k$ on  $\nu_\W(M,N)_\lin^i$ is given by 
	\[ (C^\infty(M)_{(i)}\cap \I^k)/		\big(C^\infty(M)_{(i)}\cap \I^{k+1}+C^\infty(M)_{(i+1)}\cap \I^k\big).\]
	The space of  polynomial functions of degree $k$ on  $\on{gr}(\nu(M,N))^i$ has exactly the same description. 
\end{proof}

\begin{proposition}\label{prop:iteratedweighted}
For a graded bundle $E\to N$,  with its canonical weighting along $N\subset M$, there is a canonical identification 
\[ \nu_\W(E,N)=E.\]
\end{proposition}
\begin{proof}
	This follows from the equality \eqref{eq:recovery} together with 
  $\on{gr}(C^\infty(E))=C^\infty_\pol(E)$ (by definition of the weighting). 
\end{proof}

\subsection{Homogeneous approximations}\label{subsec:homapp}
The filtration of functions on $M$ induces filtrations on all spaces of tensor fields on $M$, and we obtain corresponding homogeneous approximations. In particular, there is an 
algebra isomorphism 
\begin{equation}\label{eq:dformgrad} \on{gr}(\Omega^\bullet(M))\to \Omega^\bullet_\pol(\nu_\W(M,N)),\end{equation}
defined by the map 
taking a $p$-form $\alpha$ of filtration degree $i$ to the form $\alpha^{[i]}\in \Omega^p(\nu_\W(M,N))^{[i]}$, homogeneous of degree $i$. 
Similarly, if $X$ is a vector field of  filtration degree $i$ (i.e.,  its action on functions by Lie derivative 
raises the filtration degree  by $i$), we obtain a homogeneous approximation $X^{[i]}$ of homogeneity $i$, 
resulting in a  Lie algebra isomorphism 
\begin{equation}\label{eq:vfgrad} \on{gr}(\mf{X}(M))\to \mf{X}_\pol(\nu_\W(M,N)).\end{equation}
All of these isomorphisms are compatible with the usual operations from Cartan's calculus; for example, 
$(\d\alpha)^{[j]}=\d(\alpha^{[j]})$ and 
$(\iota_X\alpha)^{[i+j]}=\iota_{X^{[i]}}\alpha^{[j]}$ for $X\in \mf{X}(M)_{(i)}$ and $\alpha\in \Omega(M)_{(j)}$.

\section{Weighted deformation space}\label{sec:weighteddef}
The deformation space for a manifold pair $(M,N)$ is a manifold $\delta(M,N)$ with a submersion onto the real line, such that the zero fiber is the normal bundle $\nu(M,N)$, while all other fibers are copies of $M$. 
The manifold structure is such that the directions normal to $N$ are `magnified' as one approaches the zero fiber. 
Excellent introductions to deformation spaces may be found in \cite{hig:eul,hig:tan}. 
An important special case is Connes' tangent groupoid \cite{con:non}, which may be seen as the deformation space $\delta(M\times M,M)$ of the pair groupoid $M\times M$ with respect to the diagonal; here $\nu(M\times M,M)\cong TM$. 
In this section, we generalize to the weighted setting. 

\subsection{Constructions}
Let $(M,N)$ be a weighted manifold pair, with an order $r$ weighting. As in \cite{hig:alg}, and extending the algebraic definition of the weighted normal bundle $\nu_\W(M,N)$, there is an algebraic definition of the 
\emph{weighted deformation space}  as the character spectrum 
\begin{equation}\label{eq:defspace} \delta_\W(M,N)=\Hom_\alg(\on{Rees}(C^\infty(M)),\R)\end{equation}
of the \emph{Rees algebra} of the filtered algebra $C^\infty(M)$, 
\begin{equation}\label{eq:rees} \on{Rees}(C^\infty(M))=\Big\{\sum_{i\in \Z} z^{-i} f_i\colon f_i\in C^\infty(M)_{(i)}\Big\}
\subset C^\infty(M)[z,z^{-1}].
\end{equation}
For $f\in C^{\infty}(M)_{(i)}$, we denote by 
\[ \ti{f}^{[i]}\colon \delta_\W(M,N)\to \R,\]
the function defined by the Laurent polynomial $z^{-i}f$.   If $f$ has filtration degree $i$, and $g$ has filtration degree $j$, then 
\[ \wt{fg}^{[i+j]}=\ti{f}^{[i]}\ti{g}^{[j]}.\]
Let 
\[ \pi\colon  \delta_\W(M,N)\to \R\] 
be the map given by evaluation on the Laurent polynomial `$z$'; thus $\pi=\wt{1}^{[-1]}$. The map $\pi$ is surjective; its fibers are  
\begin{equation}\label{eq:fibers} \pi^{-1}(t)=\begin{cases} \nu_\W(M,N)& t=0,\\
M& t\neq 0.\end{cases}\end{equation}
This follows since an algebra morphism $\on{Rees}(C^\infty(M))\to \R$ taking $z$ to $t$ is the same as an algebra morphism from the quotient by the ideal $(z-t)\on{Rees}(C^\infty(M))$; this quotient is equal to $\on{gr}(C^\infty(M))$ if $t=0$ and to $C^\infty(M)$ if $t\neq 0$. For $f\in C^\infty(M)_{(i)}$, we obtain 
\[ \ti{f}^{[i]}\colon \delta_\W(M,N)\to \R,\ \ \ 
\ti{f}^{[i]}\big|_{\pi^{-1}(t)}=\begin{cases} f^{[i]} & t=0,\\t^{-i} f & t \neq 0.\end{cases}\]

\begin{theorem}\label{th:defspace}
The space $\delta_\W(M,N)$ has a unique $C^\infty$ structure, as a  manifold of dimension $\dim(M)+1$, in such a way that 
the functions  given by evaluation on elements of  $\on{Rees}(C^\infty(M))$ are all smooth. 	
In terms of this manifold structure,  the 
map 
\[ \pi\colon  \delta_\W(M,N)\to \R\] 
 is a surjective submersion, and the identifications \eqref{eq:fibers} are diffeomorphisms. 
\end{theorem}

\begin{proof}
The construction is a straightforward extension from the unweighted case (see \cite{hig:eul,hig:tan}). For 
weighted coordinates $x_1,\ldots,x_n$ on open subsets $U\subset M$, we shall take the functions 
\[ y_a=\wt{x}_a^{[w_a]},\ a=1,\ldots,n,\]
together with the variable `$t$' (given by the projection $\pi$) as coordinates on $\delta_\W(U,U\cap N)$. 
Denote by $\ol{U}\subset \R^n$ the image of $U$ under the coordinate map $(x_1,\ldots,x_n)\colon U\to \R^n$. 
Since 
	\[y_a=\wt{x}_a^{[w_a]}=\begin{cases}{x}_a^{[w_a]}& t=0,\\t^{-w_a} x_a&t\neq 0. \end{cases}\] 
we see that the map $(y_1,\ldots,y_n,t)\colon \delta_\W(U,U\cap N)\to \R^{n+1}$ is a bijection onto the open subset 
\[  \{(y_1,\ldots,y_n,t)\in \R^{n+1}|\ (t^{w_1}y_1,\ldots,t^{w_n}y_n)\in \ol{U}\}.\]
Clearly, a covering of $M$ by weighted coordinate charts $U$ determines a covering of $\delta_\W(M,N)$ by charts 
$\delta_\W(U,U\cap N)$. 

	
	We need to show that weighted coordinate changes on $M$ gives rise to smooth coordinate changes on the deformation space $\delta_\W(M,N)$. By considering the components of a weighted coordinate change, it is enough to show that if $U$ is the domain of a set of weighted coordinates $x_a$, and  $f\in C^\infty(U)$ is a function of filtration degree $i$, then the restriction of $\wt{f}^{[i]}$ to  $\delta_\W(U,N\cap U)$ is smooth. 
	We may  assume that $f(x_1,\ldots,x_n)=\chi(x_1,\ldots,x_n)\,x^s$ with $w\cdot s\ge i$ and $\chi\in C^\infty(U)$. 
	Expressed in deformation space coordinates, the function $\wt{f}^{[i]}$ is given for $t\neq 0$ by 
	\[
	\wt{f}^{[i]}(y_1,\ldots, y_n,t)=t^{-i} f(x_1\ldots,x_n)=t^{-i}f(t^{w_1}y_1,\ldots,t^{w_n}y_n)
	=t^{w\cdot s-i}\chi(t^{w_1}y_1,\ldots,t^{w_n}y_n) y^s.\]
	We see  that this function extends smoothly to $t=0$. The limiting function is zero if $w\cdot s>i$ (so that $f$ actually has filtration degree $i+1$), and for $w\cdot s=i$ is given by $\chi(y_1,\ldots,y_{k_0},0,\ldots,0)\ y^s$. 
	In both cases, this agrees with 
	\[ \wt{f}^{[i]}(y_1,\ldots, y_n,0)=f^{[i]}(y_1,\ldots,y_n),\]
	as required.
	
	The remaining claims are immediate from the coordinate description. 
\end{proof}

\begin{remark}
	The construction may be adapted to the analytic or holomorphic settings, by working with the sheaf of algebras $\on{Rees}(C^\infty_M)$. For $N=\pt$, a weighted deformation to the normal cone was used in Kaveh's work on toric degenerations \cite{kav:tor}.
\end{remark}

The construction of weighted deformation spaces is functorial: a morphism of weighted pairs
$\varphi\colon (M,N)\to (M',N')$ defines a smooth map of deformation spaces 
\[ \delta_\W(\varphi)\colon \delta_\W(M,N)\to \delta_{\W}(M',N'),\]
given on the open piece by $M\times \R^\times \to M'\times \R^\times, (m,t)\mapsto (\varphi(m),t)$, and on 
the zero fiber  by $\nu_\W(\varphi)\colon \nu_\W(M,N)\to \nu_{\W}(M',N')$. (In terms of the algebraic definition, it is induced by 
the morphism of Rees algebras $\on{Rees}(C^\infty(M'))\to \on{Rees}(C^\infty(M))$ given by the pullback of functions.) 
The smoothness of $\delta_\W(\varphi)$ follows since for any 
$g\in C^\infty(M')_{(i)}$,  the pullback of $\wt{g}^{[i]}$ is the smooth function $\wt{g\circ \varphi}^{[i]}$, while the 
pullback of the function `$t$' on $M'$ (given by the projection $\pi$ to $\R$) is the function `$t$' on $M$.  As a special case, we can apply this to the map of manifold pairs 
$(N,N)\to (M,N)$, where $(N,N)$ has the trivial weighting. Since $\delta(N,N)=N\times \R$, this gives a 
natural embedding 
\[ N\times \R\to \delta_\W(M,N).\]

\subsection{Basic properties}
Here are some other important aspects of the weighted deformation space. 

\begin{enumerate}
	\item If $E\to N$ is a graded bundle, with its canonical weighting, then $\delta_\W(E,N)\cong E\times \R$. 
	\item The  map 
	\[ \delta_\W(M,N)\to M, \]
	given on $M\times \R^\times $ by projection to the first factor and on $\nu_\W(M,N)$ by bundle projection to $N$, followed by inclusion, is smooth.  This follows from the local coordinate description, or since its composition  
	with any $f\in C^\infty(M)$ is the smooth function 
	$\wt{f}^{[0]}$. 
	\item 
	The space $\delta_\W(M,N)$ comes with an action $u\mapsto \wt{\kappa}_u$ of the non-zero scalars, in such a way that each 
	$\ti{f}^{[i]}$ is homogeneous of degree $i$. The action of $u\in \R^\times$ is defined  by the algebra morphism 
	$\sum_{i\in \Z} z^{-i} f_i\mapsto \sum_{i\in \Z} u^i z^{-i} f_i$ of the Rees algebra. 
	%
	The action map
	\[ \wt{\kappa}\colon \R^\times \times \delta_\W(M,N)\to \delta_\W(M,N),\ \ (u,p)\mapsto \wt{\kappa}_u(p)\]
	is smooth, since its 
	composition with  $\wt{f}^{[i]}$, for a function $f$ of filtration degree 
	$i$, is the smooth function $(u,p)\mapsto \wt{f}^{[i]}(u\cdot p)=u^i \wt{f}^{[i]}(p)$, and its composition with $\pi$ is the smooth function $(u,p)\mapsto u^{-1}\cdot \pi(p)$. 
	\item 
	The 
	vector field 
	\[ \Theta\in\mf{X}(\delta_\W(M,N))\] with flow $\tau\mapsto \wt{\kappa}_{e^\tau}$ 
	is given by $-\E$ on $\nu_\W(M,N)$ and by $t\f{\p}{\p t}$ on $M\times \R^\times$. 
\end{enumerate}

\subsection{Homogeneous interpolations}
For $f\in C^\infty(M)_{(i)}$, the function $\wt{f}^{[i]}\in C^\infty(\delta_\W(M,N))$ interpolates between the original function $f$ at $t=1$ and its homogeneous approximation $f^{[i]}$ at $t=0$. As already mentioned, it is 
homogeneous of degree $i$ for the $\R^\times$-action.  We get similar interpolations for all tensor fields of a given filtration degree.  
For example, given a vector field $X$ of  filtration degree $i\ge -r$, then the vector field 
$t^{-i} X$ on the open subset $M\times \R^\times$ extends to a vector field $\wt{X}^{[i]}$ on $\delta_\W(M,N)$. If $\alpha\in \Omega(M)$ has filtration degree $i$ then $t^{-i}\alpha$ extends to a form 
$\wt{\alpha}^{[i]}$. These extensions are homogeneous of degree $i$, and along the zero fiber $\nu_\W(M,N)$, one obtains the homogeneous approximations 
$X^{[i]}$, respectively $\alpha^{[i]}$. The Lie derivative $\L_\Theta$ acts as $-i$ on these extensions, for example, 
\[ \L_\Theta \wt{\alpha}^{[i]}=-i\wt{\alpha}^{[i]}\]
if $\alpha$ has filtration degree $i$. Once again, the extensions are compatible with the Cartan calculus. For example, in order to prove  $\L_{\wt{X}^{[i]}}\wt{\alpha}^{[j]}=(\wt{\L_X\alpha})^{[i+j]}$, 
it suffices to remark that the formula holds over the open set $M\times \R^\times$. 

\begin{example}
	Let $\pi\in\mf{X}^2(M)$ be a Poisson bivector field, that is, the Schouten bracket with itself is zero: $[\pi,\pi]=0$. 
	If $\pi$ has filtration degree $i$ for a given weighting of $(M,N)$, then we obtain a bivector field on the deformation space, 
	homogeneous of degree $i$, 
	\[ \wt{\pi}^{[i]}\in\mf{X}^2(\delta_\W(M,N)),\]
	which is given by $t^{-i} \pi$ on all non-zero fibers and by the homogeneous approximation $\pi^{[i]}$ on the zero fiber $\nu_\W(M,N)$. Both 
	$\wt{\pi}^{[i]}$ and $\pi^{[i]}$ are again Poisson structures, since Schouten bracket of the $i$-homogeneous approximations is the $2i$-homogeneous approximation of the Schouten bracket $[\pi,\pi]$, hence is zero. 
\end{example}

\subsection{Euler-like vector fields}

\begin{definition}
	A vector field $X\in\mf{X}(M)$ is called (weighted) \emph{Euler-like} for the weighted manifold pair $(M,N)$ 
	if it has filtration degree $0$, with homogeneous approximation $X^{[0]}$ the Euler vector field on the weighted normal bundle. 
\end{definition}

Recall that in graded bundle coordinates, the Euler vector field on a graded bundle with weights $w_1,\ldots,w_n$ is given by 
$\E=\sum_{a=1}^n w_a y_a\f{\p}{\p y_a}$.

\begin{definition}
	A	smooth map $\varphi\colon \nu_\W(M,N)\supset O\to M$, defined on a star-shaped open neighborhood $O$ of the zero section $N$, is called a \emph{weighted tubular neighborhood embedding} if it is a morphism of weighted manifolds, with homogeneous approximation the identity map of $\nu_\W(M,N)$. 
\end{definition}
Here we are using the natural identification $\nu_\W(E,N)=E$ (see Proposition \ref{prop:iteratedweighted})  for $E=\nu_\W(M,N)$. In analogy with the key observation in \cite{bur:spl} (see also \cite{bis:def,dan:ext,hig:eul}), we have: 
\begin{theorem}\label{th:euler}
	Every Euler-like vector field $X$ for the weighted manifold pair $(M,N)$ defines a unique maximal tubular neighborhood embedding $\varphi\colon \nu_\W(M,N)\supset O\to M$ such that $\varphi^*(X|_{\varphi(O)})=\E|_O$. If $X$ is complete, one can take $O=\nu_\W(M,N)$.
\end{theorem}

Here the uniqueness part means that any other tubular neighborhood embedding with these properties is obtained by restriction to an open subset of $O$. The proof is parallel to the unweighted case \cite{hig:eul} (cf. \cite{bis:def}): One observes that the vector field $\f{1}{t}(\wt{X}^{[0]}+\Theta)$ is well-defined;  the  diffeomorphism $\varphi$ is obtained by applying its 
flow to the zero fiber $\nu_\W(M,N)$. 

\begin{remark}
	Note that the weighting for a manifold pair  $(M,N)$ may be recovered from the corresponding Euler-like vector fields $X$. 
	Indeed, letting $s\mapsto \Phi_s^X$ be the flow of $X$ (defined on $O$, for all $s\ge 0$), the maps $\lambda_t=\Phi^X_{-\log(t)}$ extend smoothly to $t=0$, and $f$ has filtration degree $i$ if and only if $\lambda_t^*f=O(t^i)$. 
\end{remark}

Theorem \ref{th:euler} shows that weighted tubular neighborhood embeddings and weighted Euler-like vector fields
are (essentially) the same thing. This is very useful for normal form problems. For example, given a tensor field $\alpha$ 
(for example, a differential form) of filtration degree $i$, one may ask for a weighted tubular neighborhood embedding $\varphi$ taking $\alpha$ to its homogeneous approximation $\alpha^{[i]}$. Using the theorem, this amounts to constructing 
a weighted  Euler-like vector field $X$ satisfying $\L_X\alpha=i\alpha$. Indeed, the pullback of this equation under the resulting $\varphi$ reads as $\L_\E\varphi^*\alpha=i\varphi^*\alpha$, which means that $\varphi^*\alpha$ is homogeneous of degree $i$, and hence coincides with $\alpha^{[i]}$. See \cite{me:eul} for some applications of this technique, in the case $r=2$.

\section{Weighted blow-ups}\label{sec:weightedblowup}
Given a submanifold $N\subset M$, its real blow-up $\on{Bl}(M,N)$ is a manifold with boundary, obtained by replacing $N$ with the sphere bundle of its normal bundle. Real blow-ups along submanifolds are widely used in the analysis of singular spaces, pioneered in the work of Melrose \cite{mel:ati}. One can also consider projective blow-ups (using 
the projectivization of the normal bundle) avoiding the introduction of a boundary,  but often introducing 
orientability issues. In this section, we will generalize to  the weighted setting. Of course, weighted 
blow-ups are widely used in algebraic geometry, and have also appeared in  $C^\infty$-differential geometry (notably in the work of Melrose and collaborators, e.g.,  \cite{kot:gen,mel:cor}).  However, we are not aware of a published account of a 
general coordinate-free framework.  

\subsection{Blow-up}
Given any graded bundle $E\to N$, we can form its sphere bundle\[ S(E)=(E- N)/\R_{>0},\]where $\R_{>0}$ is the multiplicative group of positive real numbers. This is indeed a fiber bundle with fibers diffeomorphic to spheres: After choice of a linearization $E\cong E_\lin$, and fixing a Euclidean fiber metric on each graded summand,  we see that each $\R_{>0}$-orbit intersects the unit sphere, consisting of vectors $v\in E\cong E_\lin$ with $||v||=1$, 
in a unique point. We may interpret $S(E)$ as the set of `closed rays' $\R_{\ge 0}\cdot v,\ v\not\in N$.  
\begin{remark}By contrast, the projectivization $P(E)=(E- N)/\R^\times$ may have orbifold singularities: in fact, it is smooth if and only if all weights from the weight sequence are odd, or all are even. In terms of a linearization $E\cong E_\lin$, this follows because nonzero vectors in $E_\lin^{-i}$ have trivial stabilizer if $i$ is odd, and stabilizer $\Z_2$ is $i$ is even.  \end{remark}

Given a weighted manifold pair $(M,N)$, we wish to define a \emph{real weighted blow-up} 
\[ \on{Bl}_\W(M,N)=S(\nu_\W(M,N))\sqcup (M-N)\] 
as a manifold with boundary.  It is convenient to consider its `doubled' version $\widehat{\on{Bl}}_\W(M,N)$
as a manifold without boundary. The following description of weighted blow-ups
in terms of weighted deformation spaces extends a description in \cite[Definition 2.9]{deb:blo}.


\begin{proposition}
	The action of $\R_{>0}$ on the subset $\delta_\W(M,N)-(N\times \R)$ is free and proper. Hence, the quotient by this action, \[ \widehat{\on{Bl}}_\W(M,N)=(M-N)\sqcup S(\nu_\W(M,N))\sqcup (M-N),\]
	inherits a manifold structure, with a $\Z_2$-action interchanging the two copies of $M-N$ and   
	preserving $S(\nu_\W(M,N))$. 
\end{proposition}
\begin{proof}
	It is obvious that the action on this subset is free. 	To see that it is proper (and hence, a principal action), it suffices to construct local slices for the action. Consider local weighted coordinates $x_1,\ldots,x_n$ on $U\subset M$, and the resulting coordinates $y_1,\ldots,y_n,t$ on $\delta_\W(U,N\cap U)$. The complement of 
	$(N\cap U)\times \R$ in this coordinate chart is characterized by the condition that at least one of the $y_a$'s with $a>k_0$ is non-zero. Define a covering of this complement by invariant open subsets
	$V^+_{k_0+1},V^-_{k_0+1}\ldots,V^+_n,V^-_n$, where $V_a^\pm$ are defined by the condition $\pm y_a>0$. On $V_a^\pm$, the $\R_{>0}$-action has a slice  given by the subset where $y_a=\pm 1$. 
\end{proof}
The manifold with boundary $\on{Bl}_\W(M,N)$ is realized as the 
image of $\pi^{-1}(\R_{\ge 0})-(N\times \R_{\ge 0})$ under the quotient map. 

\begin{remark}
	The action of 	$\R_{>0}$ on  $\delta_\W(M,N)-(N\times \{0\})$ is also free. However, it is not proper, and the quotient space fails to be a manifold.
\end{remark}

\subsection{Blow-down}
Observe that the canonical map $\delta_\W(M,N)\to M$ is $\R^\times$-invariant. It hence descends to a smooth \emph{blowdown map}
\[ \widehat{\on{Bl}}_\W(M,N) \to M.\]
On each of the two copies of $M-N$, the blow-down map is the obvious inclusion, and on the exceptional divisor  
it is the bundle projection $S(\nu_\W(M,N))\to N$ followed by inclusion $N\to M$. 

\begin{proposition}
	The blow-down map is a morphism of weighted manifold pairs, 
	\[ 	 \big(\widehat{\on{Bl}}_\W(M,N),S(\nu_\W(M,N))\big)
	\to (M,N)\]
	using the trivial weighting for the pair  $\big(\widehat{\on{Bl}}_\W(M,N),S(\nu_\W(M,N))\big)$. 
	Every Euler-like vector field $X\in\mf{X}(M)$, with respect to the given weighting of $(M,N)$, lifts to an Euler-like vector field on $\widehat{\on{Bl}}_\W(M,N)$, with respect to the trivial weighting . 
\end{proposition}	
\begin{proof}
	It is enough to prove the second claim, since weightings are uniquely determined by Euler-like vector fields. Using the tubular neighborhood embedding defined by $X$, we may assume that $M=E$ is a graded bundle, with $X=\E$ \emph{the} Euler vector field of $E$. Hence, in graded bundle coordinates $x_1,\ldots,x_n$, 
	\[ X=\sum_a w_a x_a\f{\p}{\p x_a}.\]
	The vector field $\wt{X}^{[0]}$, written in coordinates $y_1,\ldots,y_n,t$, is given by a similar expression 
	$\wt{X}^{[0]}=\sum_a w_a y_a\f{\p}{\p y_a}$. We want to express this vector field in a chart $V_c^\pm$, 
	with coordinates 
	\[ z_a=y_a y_c^{-w_a/w_c},\,a\neq c,\ \ \ z_c=t y_c^{1/w_c},\ \ t.\]
	(Note that in the new coordinates, the action of $\R_{>0}$ on the $z_a$'s becomes trivial.)  An elementary calculation gives $\wt{X}^{[0]}=z_c \f{\p}{\p z_c}$, confirming the claim.  
\end{proof}	

\subsection{Further remarks}
Some further comments on weighted blow-ups: 
\begin{enumerate}
	\item Tensor fields on $M$ lift to the weighted blow-up if and only if they are of filtration degree $0$. 
	For example, if $X\in\mf{X}(M)$ has filtration degree $0$, then the vector field $\wt{X}^{[0]}$ on 
	$\delta(M,N)-(N\times \R)$ is $\R^\times$-invariant, and so descends. Similarly, if $\alpha\in \Omega(M)_{(0)}$ then $\wt{\alpha}^{[0]}$ on 
	$\delta(M,N)-(N\times \R)$ is $\R^\times$-basic, and so it descends.  
	\item 
	We defined the smooth structure on the weighted blow-up in terms of the deformation space. 	One can also proceed in the 
	other direction, viewing $\delta_\W(M,N)$ as a submanifold of a weighted blow-up.  In fact, there is a canonical decomposition 
	\[ \wh{\on{Bl}}_\W(M\times\R,N\times \{0\})=\delta_\W(M,N)\sqcup \wh{\on{Bl}}_\W(M,N)\sqcup \delta_\W(M,N),\]
	with $\wh{\on{Bl}}_\W(M,N)$ embedded as a hypersurface and the two copies of the deformation space embedded as open subsets.
	\item Suppose   $\varphi\colon (M,N)\to (M',N')$ is a morphism of weighted manifold pairs. In order to lift to the blow-up manifolds, it is necessary to remove a certain closed subset: we obtain a smooth lift $\widehat{\on{Bl}}^0_\W(M,N)\to \widehat{\on{Bl}}_{\W}(M',N')$ 
	where 
	\[ \widehat{\on{Bl}}^0_\W(M,N)=\big(\delta_\W(M,N)-\delta_\W(\varphi)^{-1}(N'\times \R)\big)/\R_{>0}
	\subset  \widehat{\on{Bl}}_\W(M,N)
	.\]
\item 
Given a linear weighting of a vector bundle $V\to M$ along a subbundle $W\to N$,  
our constructions give vector bundles
\[ \nu_\W(V,W)\to \nu_\W(M,N),\ \ \delta_\W(V,W)\to \delta_\W(M,N),\ \ \wh{\on{Bl}}^0_\W(V,W)\to \wh{\on{Bl}}_\W(M,N),\]
where $ \wh{\on{Bl}}^0_\W(V,W)\subset \wh{\on{Bl}}_\W(V,W)$ is the quotient of $\delta_\W(V,W)|_{\delta_\W(M,N)-N\times\R}$ by $\R^\times$. 
\end{enumerate}

\section{Weightings as subbundles of the $r$-th tangent bundle}\label{sec:weightingstrm}
For $r>2$, the filtration  of the normal bundle $\nu(M,N)$   
does \emph{not} suffice to describe a weighting along $N\subset M$. Intuitively, the subbundles $F_{-i}$ 
in \eqref{eq:normalbundlefiltration}
only give first-order information, whereas the weighting requires higher order information.  This motivates working with jet spaces.  We will see that order $r$ weightings are in 1-1 correspondence with certain subbundles $Q$ of the $r$-th tangent bundle $T_rM$, and the weighted normal bundle is a quotient of that subbundle. For a trivial weighting, we have $Q=TM|_N$, and we recover the 
usual description of the normal bundle as a quotient $TM|_N/TN$. 

\subsection{Higher tangent bundles}
We will need additional background material on higher tangent bundles $T_rM$. 
\subsubsection{Algebraic definition}
In Example \ref{ex:gradedbundles}\eqref{it:e} we recalled the definition of $T_rM$  as the space $J^r_0(\R,M)$
of $r$-jets of curves. Similar to the description of tangent vectors as derivations $v\colon C^\infty(M)\to \R$, there is  a more algebraic approach to $T_rM$, which will be more convenient for our purposes. See  \cite[Chapter VIII]{kol:nat} for detailed discussions and proofs. 
Consider the truncated polynomial algebra
\begin{equation}\label{eq:aar}
\AA_r=\R[\epsilon]/(\epsilon^{r+1}),
\end{equation}
equipped with the grading where $\epsilon^i$ has degree $-i$. 
The $r$-th tangent bundle is the space of algebra morphisms, 
\begin{equation} \label{eq:trm} T_rM=\Hom_\alg(C^\infty(M),\AA_r).\end{equation}
Elements of $T_rM$ may be written $\su=\sum_{i=0}^r \su_i \epsilon^i$, with linear maps $\su_i\colon C^\infty(M)\to \R$.  Here $\su_0$ is an algebra morphism (hence is given by evaluation at some point $m\in M$), $\su_1$ is a derivation with respect to $\su_0$ (hence is given by a tangent vector based at $m$), and so on. 
The identification of the jet space $J_0^r(\R,M)$ with the space $\Hom_\alg(C^\infty(M),\AA_r)$ is  the map taking  
$\sj^r_0(\gamma)$ to the algebra morphism 
\[
C^\infty(M)\to \AA_r,\ \ f\mapsto \sum_{i=0}^r \f{1}{i!}\f{d^i}{d t^i}\Big|_{t=0}f(\gamma(t))\epsilon^i.
\]
Put differently, this is the algebra morphism 
$f\mapsto \sj^r_0(f\circ \gamma)$, where we identify $J_0^r(\R,\R)\cong \AA_r$.

\begin{remark}
The same terminology `$r$-th tangent bundle'
if often used for  the bundle 
\[ T^rM=\Hom_\alg(C^\infty(M),\AA_1\otimes \cdots \otimes \AA_1).\] We shall call the latter the \emph{$r$-fold tangent bundle} since it is isomorphic to the iteration $TT\cdots TM$.  The symmetric group $S_r$ acts on $\AA_1\otimes \cdots \otimes \AA_1$ by permutation of the factors. Its fixed point set is $\AA_r$, embedded by the map taking the generator $\epsilon$ of 
$\AA_r$ to the sum $\sum_i \epsilon_i$ of the generators for various copies of $\AA_1$. This gives an embedding $T_rM\to T^rM$ as the `symmetric part' of the $r$-fold tangent bundle. 
\end{remark}

In the algebraic picture, the $(\R,\cdot)$-action on $T_rM$ is given by $\kappa_t\colon \sum\su_j \epsilon^j\mapsto 
\sum_j\su_j t^j  \epsilon^j$
(cf.~ Example \ref{ex:gradedbundles}\eqref{it:b}). It extends further to a monoid action of 
\begin{equation}\label{eq:Lambdar}
\Lambda_r=\Hom_\alg(\AA_r,\AA_r),\end{equation}
by composition, $\Psi.\mathsf{u}=\Psi\circ \mathsf{u}$.  (See \cite{joz:monoid} for a detailed discussion of the $\Lambda_2$-action on $T_2M$.) As a vector space, $\Lambda_r\cong \epsilon \AA_r=\AA_r^-$, since 
any element $\Psi\in \Lambda_r$ is uniquely determined by the image of $\epsilon$, and since the degree $0$ part of $\Psi(\epsilon)$ must vanish due to $\Psi(\epsilon)^{r+1}=\Psi(\epsilon^{r+1})=0$.  
In the jet picture, $\Lambda_r=J^r_{0,0}(\R,\R)^{\on{op}}$ acting on $J^r_0(\R,M)$ by composition. 

The quotient maps 
\[ p^r_s\colon T_rM\to T_sM\]
for $r\ge s$ are induced by  the algebra morphisms $\AA_r\to \AA_s$, and in particular $p^r_0$ is the base projection
$T_rM\to M$. 
The $r$-th tangent lift of a smooth map $F\colon M\to M'$ is the morphism of graded bundles  
$T_rF\colon T_rM\to T_rM'$, taking $\mathsf{u}\in T_rM$ to $\mathsf{u}\circ F^*$, where $F^*$ is the pullback of functions. 
In particular,  diffeomorphisms of $M$ act on $T_rM$ by graded bundle automorphisms.

\subsubsection{Lift of functions}\label{subsubsec:functions}
Every function $f\in C^\infty(M)$ determines a sequence of functions $ f^{(i)}\in C^\infty(T_rM)^{[i]}$, homogeneous of degree $i$,  
taking $\su=\sum_{j=0}^r \su_j \epsilon^j\in T_rM$ to 
\begin{equation}\label{eq:fidef} f^{(i)}(\su)=\su_i(f).\end{equation}
In the jet picture, the map takes $\su=\sj^r_0(\gamma)$ to $\f{1}{i!}\f{d^i}{d t^i}\big|_{t=0}f(\gamma(t))$.
In particular, $f^{(0)}$ is  the pullback $(p^r_0)^* f$, while $f^{(1)}$ is the pullback under $p^r_1$ of 
the exterior differential $\d f$, regarded as a function on $TM$. Under multiplication of functions, one has the \emph{product rule} 
\begin{equation}\label{eq:productrule}
(fg)^{(i)}=\sum_{j=0}^i f^{(j)}g^{(i-j)}.\end{equation}
If $x_a$ are local coordinates 
on $U\subset M$, then the functions 
\[ x_a^{(i)},\ \ \ a=1,\ldots,n,\ i=0,\ldots,r\] serve as graded bundle coordinates on $T_rU\subset T_rM$. The function $f^{(i)}$  is a polynomial in 
the variables $x_a^{(j)}$ with $0<j\le i$. In low degrees, 
\begin{align*}
f^{(1)}&=\sum_a \f{\p f}{\p x_a} x_a^{(1)},\\
f^{(2)}&=\sum_a \f{\p f}{\p x_a} x_a^{(2)}
+\f{1}{2}\sum_{a_1a_2} \f{\p^2 f}{\p x_{a_1} \p x_{a_2}} x_{a_1}^{(1)} x_{a_2}^{(1)},\\
f^{(3)}&=\sum_a \f{\p f}{\p x_a} x_a^{(3)}
+\sum_{a_1a_2} \f{\p^2 f}{\p x_{a_1} \p x_{a_2}} x_{a_1}^{(2)} x_{a_2}^{(1)}
+\f{1}{6}\sum_{a_1a_2 a_3} \f{\p^3 f}{\p x_{a_1} \p x_{a_2} \p x_{a_3}} x_{a_1}^{(1)} x_{a_2}^{(1)}  x_{a_3}^{(1)} ;
\end{align*} 	
these expressions are found as the Taylor coefficients of $t\mapsto f(x(t))$. As one may see from these formulas, the
 $C^\infty(M)$-module of homogeneous polynomials spanned by functions of the form $f^{(i)}$ is a proper submodule of $C^\infty(T_rM)^{[i]}$ if $n\ge 2$. For example, $x_1^{(1)}x_2^{(2)}$ is not in the span, since the formula for $f^{(3)}$ only gives symmetric combinations such as $x_1^{(1)}x_2^{(2)}+x_2^{(1)}x_1^{(2)}$.

\subsubsection{Lifts of vector fields}\label{subsubsec:action}
As observed by A.~ Morimoto \cite{mor:lif}, every vector field  $X\in \mf{X}(M)$  gives rise to vector fields
$ X^{(-i)}\in \mf{X}(T_rM)^{[-i]}$, for $ i=0,\ldots,r$, on the r-th tangent bundle. These are uniquely determined by their actions on 
 lifts of functions $f\in C^\infty(M)$, as follows: 
\begin{equation}\label{eq:liftdef}
X^{(-i)}f^{(\ell)}= (X f)^{(\ell-i)}\end{equation}
The vector field $X^{(0)}$ is the  \emph{r-th tangent lift}, while $X^{(-i)}$ for $0<i\le r$ are the \emph{vertical lifts}. 
In local coordinates, 
\begin{equation}\label{eq:xicoord} 
 \big(\sum_{a=1}^n f_a \f{\p}{\p x_a}\big)^{(-i)}=\sum_{k=i}^r  
\sum_{a=1}^n f_a^{(k-i)}\f{\p}{\p x_a^{(k)}}.\end{equation}
The product rule $(fX)^{(-i)}=\sum_{k=i}^r f^{(k-i)} X^{(-k)}$  shows that the maps  $X\mapsto X^{(-i)}|_M$ are $C^\infty(M)$-linear, since $f^{(j)}|_M=0$ for $j>0$.  This 
realizes the isomorphism of vector bundles 
\begin{equation}\label{eq:trm}
	T(T_rM)|_M=TM\oplus (T_rM)_\lin=TM^{\oplus(r+1)}\end{equation}
The lifts satisfy bracket relations 
\begin{equation}\label{eq:bracketXY} [X^{(-i)},Y^{(-j)}]=[X,Y]^{(-i-j)},\end{equation}
where the right hand side is zero if $i+j>r$.

\subsection{The graded subbundle defined by a weighting}
The following result shows that an order $r$ weighting  along a submanifold $N\subset M$ determines a graded subbundle of the $r$-th tangent bundle, and may be recovered from that subbundle. Let 
\begin{equation}\label{eq:defofq} Q=\{\su\in T_rM|\ 
\,\forall f\in C^\infty(M)_{(i)},\ j<i\le r\Rightarrow
f^{(j)}(\su)=0\}.\end{equation}
Denote by   $\ti{F}_{-i}\subset TM|_N$ the pre-images of the subbundles $F_{-i}\subset \nu(M,N)$ from Proposition \ref{prop:f}.
 
 \begin{theorem}[Weightings as subbundles of $T_rM$]\label{th:q}
 \begin{enumerate}
 	\item The subset $Q\subset T_rM$ is a graded subbundle $Q\to N$ of dimension  $k_0+\ldots+k_r$. In fact, $Q$ is $\Lambda_r$-invariant. 
 	\item The linear approximation  is 
 	\[ Q_\lin=\ti{F}_{-1}\oplus \cdots \oplus \ti{F}_{-r},\]
 	as a graded subbundle of $(T_rM)_\lin=TM^{\oplus r}$. In particular, $Q_\lin^{-r}=TM|_N$. 
 	\item The weighting is recovered from $Q$ as follows:
 	\[ 	C^\infty(M)_{(i)}=\{f\in C^\infty(M)|\ \forall j< i\colon f^{(j)}	\mbox{ vanishes on}\ Q\}\]
 	for $i\le r$. 
 \end{enumerate}
 \end{theorem}
Note that for  the case of a trivial weighting ($r=1$), the graded subbundle $Q\subset T_1M$ is the restriction $TM|_N\subset TM$. 

\medskip
 The proof of Theorem \ref{th:q} is based on the following local coordinate description of $Q$:
 %

 \begin{lemma}\label{lem:qincoord} In weighted coordinates $x_a\in C^\infty(U)$, the set $Q\cap T_rU$ is 
 	given by the system of equations
 	\begin{align*}
 	x_a^{(0)}=0& \ \ \mbox{ for } w_a>0,\\ 
 	x_a^{(1)}=0 &\ \ \mbox{ for } w_a>1,\\ 
 	\ldots &\\
 	x_a^{(r-1)}=0&\ \ \mbox{ for } w_a>r-1. 
 	\end{align*}
 	Hence, the collection of all $x_a^{(j)}$ with $w_a\le j$ serve as  coordinates on $Q\cap T_rU$.
 \end{lemma}
 \begin{proof}
 	Let $Q'\subset T_rU$ be the subset given by this system of equations.  
 	We claim $Q\cap T_rU=Q'$. For the inclusion $\subset$, note that if $w_a>j$, then $x_a$ has filtration degree 
 	$j+1$, and so $x_a^{(j)}$ must vanish on $Q\cap T_rU$.  For the opposite inclusion $\supset$,  	
 	consider the monomials $x^s=x_1^{s_1}\cdots x_n^{s_n}$ for multi-indices $s=(s_1,\ldots,s_n)$. 
 	By induction on $\sum_a s_a$, and using the product rule \eqref{eq:productrule},
 	 the function $(x^s)^{(j)}$ vanishes on $Q'$ for $\sum_a s_a w_a>j$. 
 	 By definition, the space of functions 
 	of filtration degree $j+1$ is spanned by products $x^s g$ with $\sum_a s_a w_a>j$; hence, using the product rule again, it follows that $f^{(j)}$ vanishes on $Q'$. 
 \end{proof}
 
 \begin{proof}[Proof of Theorem \ref{th:q}]
 	\begin{enumerate}
 		\item The local coordinate description shows that $Q$ is a graded subbundle, of the stated  dimension. Let $\Psi\in \Lambda_r, \su\in Q$. To show $\Psi \su\in Q$, let 
 		$f\in C^\infty(M)_{(i)}$ with $i\le r$ be given. We have, using  \eqref{eq:fidef}, 
 		\[ \sum_{j=0}^r f^{(j)}(\Psi\su)\epsilon^j=(\Psi\su)(f)=\Psi(\su(f))=
 		\sum_{j=i}^r f^{(j)}(\su)(\Psi(\epsilon))^j. \]
 		On the right hand side of this equality, since $\Psi(\epsilon)\in \epsilon\AA_r$, 
 		the coefficients of $\epsilon^j$ with $j<i$ are zero. Hence, the same is true for the left hand side, proving that 
 		$j<i\Rightarrow f^{(j)}(\Psi\su)=0$. We conclude $\Psi\su\in Q$. 
 		\item 
 		We shall use the description of $Q_\lin$ in terms of $TQ|_N=TN\oplus Q_\lin$, where the second summand is embedded as the tangent space to the fibers of $Q\to N$. Similarly $T(T_rM)|_M=TM\oplus (T_rM)_\lin$. In terms of  weighted coordinates over $U\subset M$, $T(T_rM)_\lin|_U$ is spanned by all $\f{\p}{\p x_a^{(j)}}$ with $0\le j\le r$ and $a=1,\ldots,n$, while 
 		\[ T(Q\cap T_rU)=\on{span}\Big\{\f{\p}{\p x_a^{(j)}}\big|\ w_a\le j \Big\}.\]
 		The identification $TM\cong T(T_rM)^{-j}$ from \eqref{eq:trm} is given by $X\mapsto X^{(-j)}$. Since
 		 \begin{equation}\label{eq:coordlift} \f{\p}{\p x_a^{(j)}}=\big(\f{\p}{\p x_a}\big)^{(-j)},\end{equation}	
 		it follows that the  pre-image of $Q_\lin^{-j}|_{N\cap U}$ in $TM|_{N\cap U}$ is spanned by all 
 		$\f{\p}{\p x_a}$ with $w_a\le j$. But this is exactly the subbundle $\ti{F}_{-j}$. 	(See the proof of Proposition 
 		\ref{prop:f}.) 
 		\item 
 		Suppose  that $f^{(j)}$ vanishes on $Q$, for all 
 		$j< i$.  We want to show that $f$ has filtration degree $i$. 
 		Using weighted coordinates, and using a Taylor expansion on the coordinates $x_a$ with $w_a>0$, this is equivalent to showing that all derivatives 
 		\[\f{\p^p f}{\p x_{a_1}\cdots \p x_{a_p}}\] 
 		with $w_{a_\nu}>0$ and $w_{a_1}+\cdots +w_{a_p}<i$ vanish along $N$. Let $j_1=w_{a_1},\ldots,j_p=w_{a_p}$, and put $j=j_1+\ldots+j_p$. Then the derivative above is  (up to a positive multiple) the coefficient of 
 		\[ 	x_{a_1}^{(j_1)}\cdots x_{a_p}^{(j_p)}\]
 		in $f^{(j)}$.
 		Since the $x_{a}^{(k)}$ with $k\le w_a$ are the coordinates on $Q$ the coefficient must be zero, as claimed. \qedhere
 	\end{enumerate}
\end{proof}

Suppose $(M,N)$ and $(M',N')$ are manifold pairs, with weightings of order $r,r'$. Raising one of the orders if needed, we may assume $r=r'$.  As a simple consequence of Theorem \ref{th:q},  we see that  a map $\varphi\colon M\to M'$ is a morphism of weighted manifolds if and only if the map $T_r\varphi\colon T_rM\to T_rM'$ satisfies 
$(T_r\varphi)(Q)\subset Q'$.

\subsection{The weighted normal bundle as a subquotient}
Let $Q\subset T_rM$ be the graded subbundle associated to the weighting. We shall give a description of the weighted normal bundle as a quotient of $Q$ under an equivalence relation.  
	
	\begin{theorem}\label{th:quotientpicture}
		The weighted normal bundle for an order $r$ weighting along $N$ 
		is canonically the quotient 
		\[ \nu_\W(M,N)=Q/\sim\]
		under the equivalence relation, 		
		\begin{equation}\label{eq:globalcond} q_1\sim q_2 \Leftrightarrow \ f^{(i)}(q_1)=f^{(i)}(q_2)\ \mbox{ for all }\
		f\in C^\infty(M)_{(i)},\ 0\le i\le r.
		\end{equation}
	\end{theorem}
		\begin{proof}
		By definition of $Q$, if $f\in C^\infty(M)$ has filtration degree $i+1$, then its $i$-th lift $f^{(i)}$ vanishes on $Q$. 
		Furthermore, if $f$ has filtration degree $i$ and $g$ has filtration degree $j$, then 
		\[ (fg)^{(i+j)}|_Q=f^{(i)}|_Q\ g^{(j)}|_Q\]
		by the product rule  \eqref{eq:productrule}, since all the other terms vanish on $Q$. This defines a map of graded algebras 
	\[\on{gr}(C^\infty(M))\to C^\infty_\pol(Q),\]
	and dually (applying $\Hom_\alg(\cdot,\R)$)
	 a morphism of graded bundles $Q\to \nu_\W(M,N)$. We claim that the fibers of this map are exactly the fibers of the equivalence relation \eqref{eq:globalcond}.  This will be evident from a coordinate description. Recall that in local weighted coordinates, 
	 $Q$ is given by the equations $x_a^{(i)}=0$ for $i<w_a$. Hence, $x_a^{(i)}$ for $i\ge w_a$ serve as coordinates on $Q$.
	 The following lemma shows that the quotient map simply forgets the coordinates with $i>w_a$, and the remaining coordinates 
	 $x_a^{(w_a)}$ descend to the standard coordinates on $\nu_\W(M,N)$. 
	\end{proof}
	\begin{lemma}Let $x_1,\ldots,x_n$ be weighted 
		coordinates on $U\subset M$.
		For $q_1,q_2\in Q\cap T_rU$, 
		\begin{equation}\label{eq:localcond}
		q_1\sim q_2\ \Leftrightarrow\ x_a^{(w_a)}(q_1)=x_a^{(w_a)}(q_2)\ \ \mbox{ for}\ a=1,\ldots,n.
		\end{equation}
	\end{lemma}
	\begin{proof}	
		The direction `$\Rightarrow$' follows from \eqref{eq:globalcond} by putting 
		$f=x_a$. Suppose conversely that $q_1,q_2\in Q$ satisfy the condition on the right hand side of \eqref{eq:localcond}. 
		To prove $q_1\sim q_2$, we want to show that $f^{(i)}(q_1)=f^{(i)}(q_2)$ for all 
		 $f\in C^\infty(U)_{(i)}$ with $0\le i\le r$.
		By the coordinate description of the filtration, 
		it suffices to consider functions of the form $f=\chi\,x_{a_1}\cdots x_{a_p}$  with $w_{a_1}+\ldots+w_{a_p}\ge i$ and $\chi\in C^\infty(U)$. 
		By the product rule \eqref{eq:productrule}, 
\[ f^{(i)}=\sum\chi^{(j)}	x_{a_1}^{(j_1)}\cdots x_{a_p}^{(j_p)}\]
where the sum is over $j_1,\ldots,j_p,j\ge 0$ with $j_1+\ldots+j_p+j=i$. Such a term vanishes on $Q$ unless 
$j_\nu\ge w_{a_\nu}$ for all $\nu$. So, only terms with $j_\nu= w_{a_\nu}$ will remain after restriction: 
\[ f^{(i)}|_Q=(\chi^{(0)}\,x_{a_1}^{(w_{a_1})}\cdots x_{a_p}^{(w_{a_p})})|_Q.\]
By assumption, the terms $ x_{a_j}^{(w_{a_j})}$ take on the same values at $q_1,q_2$. On the other hand, 
 $\chi^{(0)}$ is the pullback of $\chi$, and \eqref{eq:localcond} for $i=0$ shows that $q_1,q_2$ have the same base point in $Q_0=N$. Hence $\chi^{(0)}$ also 
		takes on the same value at $q_1,q_2$. 	By \eqref{eq:globalcond}, we conclude   $q_1\sim q_2$. 
	\end{proof}
	
\begin{remark}
		For a trivial weighting $r=1$, we have $Q=TM|_N\subset TM$. Here, the $i=0$ part of \eqref{eq:globalcond} 
		says that the tangent vectors $q_1,q_2$ have the same base point in $N$, and the $i=1$ part says that they agree on the vanishing ideal of 	$N$, i.e., their difference is tangent to $N$. This correctly gives the equivalence relation for $\nu(M,N)$ as a quotient $TM|_N/TN$. 
\end{remark}

	\subsection{Definition of weighted normal bundle in terms of filtration of vector fields}\label{subsec:gradedliealgebroid}
	Given a weighting on $M$ along $N$, consider the resulting filtration on the Lie algebra of vector fields. The
	associated graded algebra is realized as polynomial vector fields on the weighted normal bundle, by the isomorphism  \eqref{eq:vfgrad}. We see that the sheaf $\on{gr}(\mf{X}(M))^-$, given as the sum of components of strictly negative degree, 
	 is the sheaf of   sections of a 	negatively graded Lie algebra bundle $\k\to N$. 
    Weighted coordinates $x_a$  over $U\subset M$ define a local frame for $\k$, given by all
    \[ x^s \f{\p}{\p x_a}\Big|_{N\cap U}\]
    with $w\cdot s-w_a<0$ and $a>k_0$, and 
    where $s$ ranges over multi-indices with $s_b=0$ for $w_b=0$.
	Similarly, the subsheaf 	
	$\on{gr}(\I\mf{X}_M)^-$ (where $\I\subset C^\infty_M$ is the vanishing ideal of $N$) is the sheaf of sections of a sub-Lie algebra bundle $\mf{l}\to N$.  A local frame for $\mf{l}$ is obtained from that for $\k$ by the additional condition $s\neq 0$. Let $L\subset K$ be the 
	bundles, over $N$,  of nilpotent Lie groups 
	exponentiating these bundles of Lie algebras. 
	
	The isomorphism \eqref{eq:vfgrad} identifies sections of $\k$ with polynomial vector fields on $\nu_\W(M,N)$ of strictly negative degree; the sections of $\mf{l}$ are polynomial vector fields which furthermore vanish along $N$. Note that 
	$\mf{X}_\pol(\nu_\W(M,N))^-$ acts fiberwise transitively on $\nu_\W(M,N)$, with stabilizer along $N$ the polynomial 
	vector fields vanishing along $N$. It follows that $K$ acts fiberwise transitively on $\nu_\W(M,N)$, with stabilizer 
	along $N$ the subbundle $L$. 
	To summarize:  
	\begin{proposition}\label{prop:klaresections}
		Given a weighted manifold pair $(M,N)$, the  $C^\infty_N$-modules \[ \on{gr}(\I\mf{X}_M)^-\subset \on{gr}(\mf{X}_M)^-\]
		are the sheaves of sections of negatively graded  Lie algebra
		bundles  $\mf{l}\subset\k$ over $N$. Letting $L\subset K$ be the corresponding bundles of nilpotent Lie groups, we have that 
		\[ \nu_\W(M,N)=K/L.\]
	\end{proposition}
In the case of a trivial weighting ($r=1$), the filtration on vector fields starts with $\mf{X}_{M,(-1)}=\mf{X}_M$, while 
$\mf{X}_{M,(0)}$ are vector fields tangent to $N$. Hence, in this case we directly have $\k=TM|_N/TN$ (with zero bracket)
	and $\mf{l}=0$.

\section{A necessary and sufficient condition}
\subsection{Some examples}
In the previous section, we saw that a weighting of order $r$ determines a graded subbundle $Q\subset T_rM$, and is uniquely determined by $Q$. However, we did not specify which graded subbundles correspond to weightings. 
In particular,  $\Lambda_r$-invariance does not suffice:
\begin{example}
	The graded subbundle $Q\subset T_2\R^2$ along $N=\R$ defined by $x_2=0,\ x_1^{(1)}=x_2^{(2)}=0$ is $\Lambda_2$-invariant, and satisfies 
	$Q_\lin^{-2}=TM|_N$,  
	but it does not correspond to a weighting since $Q_\lin^{-1}$ is the zero bundle (and so does not contain $TN$).
\end{example}
We might add the condition that the summands of $Q_\lin\subset (T_rM)_\lin=TM^{\oplus r}$ define a filtration of $TM|_N$, 
with $Q_\lin^{-r}=TM|_N$, but this still is not enough:   
\begin{example}\label{ex:counterexample}
	The graded subbundle $Q\subset T_4\R^3$ along $N=\{0\}$, given by the equations 
	\[ x_1=x_2=x_3=0,\ x_3^{(1)}=0,\ x_3^{(2)}=0,\ x_3^{(3)}=x_1^{(1)}x_2^{(2)}-x_1^{(2)}x_2^{(1)},\]
	is $\Lambda_4$-invariant, and defines a filtration of $TM|_N$, 
	but it does not correspond to a weighting. Indeed, if $Q$ were the graded subbundle defined by a weighting, 
	then the prescription from Theorem \ref{th:q} would show that $x_1,x_2$ have filtration degree $1$ (but not $2$) while $x_3$ has filtration degree $3$ (but not $4$).  However, 
	the graded subbundle  defined by this weighting is given by the equations $x_1=x_2=x_3=0,\ x_3^{(1)}=0,\ x_3^{(2)}=0$, and so is strictly larger than $Q$.  
\end{example}

\subsection{Infinitesimal transitivity}
The key properties of graded subbundles corresponding to weightings come from the following result. 

\begin{proposition} \label{prop:gq}
	Let $Q\subset T_rM$ be the graded subbundle defined by a weighting of order $r$. 
	\begin{enumerate}
		\item The lift $X^{(-i)}$ of a vector field $X\in\mf{X}(M)$ is tangent to $Q$ if and only if $X$ has filtration degree $-i$. 
		\item 
		Vector fields of the form $X^{(-i)}$ with $0\le i\le r$ and $X\in\mf{X}(M)_{(-i)}$ span the tangent bundle of  
		$Q$ everywhere. 
	\end{enumerate}
\end{proposition}
\begin{proof}
	We use  local weighted coordinates $x_a$ on $U\subset M$. A vector field   $X=\sum f_a \f{\p}{\p x_a}\in\mf{X}(U)$ has filtration degree $-i$  if and only if each $f_a$ has filtration degree $w_a-i$. 
	On the other hand, $X^{(-i)}$ is given in coordinates $x_a^{(k)}$ by  \eqref{eq:xicoord}. Since 
	\[ T(Q\cap T_rU)=\on{span}\big\{\f{\p}{\p x_a^{(j)}}\big|\ w_a\le j \big\},\]
	we see that $X^{(-i)}$ is tangent to $Q$ if and only if the coefficients  $f_a^{(j-i)}$ with 
	$j<w_a$ vanish on $Q$. This is equivalent to saying that each $f_a$ has filtration degree $w_a-i$, proving  (a). 
	Part (b) is also clear from local coordinates, using \eqref{eq:coordlift}. 
\end{proof}
Recall that there is a natural vector bundle action of $TM\to M$  on  $T_rM\to M$. In terms of the description 
$T_rM=\Hom_\alg(C^\infty(M),\AA_r)$, the action of a tangent vector $v\in TM|_m$ (regarded as a linear map $C^\infty(M)\to \R$ with the usual derivation property) 
on $\su\in T_rM|_m$ 
is given by 
\[ \su\mapsto \su+v\epsilon^r.\]
Note that the action of $X\in \mf{X}(M)=\Gamma(TM)$ is the time-1 flow of the $r$-th vertical lift $X^{(-r)}$.
If $Q\subset T_rM$ corresponds to a weighting, we saw that all such lifts are tangent to $Q$ (since 
$\mf{X}(M)=\mf{X}(M)_{(-r)}$). That is, $Q$ is invariant under the $TM$-action on $T_rM$. 


We shall also need the following 
natural action of the algebra $\AA_r$ on the tangent spaces of $T_rM$. 
This action was observed by Koszul; for detailed discussions see the articles of  A.~Morimoto \cite{mor:pro}
and  Okassa \cite{oka:pro}. One possible description uses the identification 
\[ T(T_rM)=\Hom_\alg(C^\infty(M),\AA_r\otimes \AA_1).\]
The action of  $x\in \AA_r$ on this space is by composition with the 
algebra endomorphism of $\AA_r\otimes \AA_1$, given on generators by 
$\epsilon\otimes 1\mapsto \epsilon\otimes 1,\ 1\otimes \epsilon \mapsto x\otimes \epsilon$.  
For a more explicit description, let us compute the action of $\epsilon$ on $X^{(-i)}|_\su\in T(T_rM)|_\su$, for 
$X\in \mf{X}(M)$ and $\su\in T_rM$. By definition (see \eqref{eq:liftdef}), $X^{(-i)}|_\su$ is the 
algebra morphism 
\[ f\mapsto \sum_j f^{(j)}(\su)(\epsilon^j\otimes 1)+\sum_j (Xf)^{(j-i)}(\su)(\epsilon^j\otimes \epsilon).\]
Hence $\epsilon\cdot X^{(-i)}|_\su$ is the algebra morphism 
\[ f\mapsto \sum_j f^{(j)}(\su)(\epsilon^j\otimes 1)+\sum_j (Xf)^{(j-i)}(\su)(\epsilon^{j+1}\otimes \epsilon).\]
Shifting indices in the second sum, this is recognized as $X^{(-i-1)}|_\su$. 
This shows that  
\begin{equation}\label{eq:araction} \epsilon\cdot X^{(-i)}=X^{(-i-1)}.\end{equation}


Consider again a graded subbundle $Q\subset T_rM$ defined by a weighting. Since $ \mf{X}(M)_{(-i)}\subset \mf{X}(M)_{(-i-1)}$,
the proposition shows that if $X^{(-i)}$ is tangent to $Q$ then so is \eqref{eq:araction}. Since the tangent bundle is spanned by such lifts, we conclude that $TQ\subset T(T_rM)|_Q$ is invariant under the action of $\AA_r$. 

\subsection{Necessary and sufficient conditions}
We are now in position to formulate a necessary and sufficient condition characterizing  weightings:
\begin{theorem}\label{th:intrinsic}
A graded subbundle  $Q\subset T_rM$ 
corresponds to a weighting if and only if the following three conditions are satisfied. 
\begin{enumerate}
	\item[(i)] \label{it:i}  $Q$ is $TM$-invariant, 
	\item[(ii)] \label{it:ii}  $TQ$ is $\AA_r$-invariant,
	\item[(iiii)] 	\label{it:iii} $TQ$ is spanned by lifts $X^{(-i)}$ of vector fields $X$ such that $X^{(-i)}$ is tangent to $Q$.  
\end{enumerate}
\end{theorem}


The necessity of these conditions was proved above.  For the other direction, suppose that a graded subbundle 
$Q\subset T_rM$ satisfying the conditions of the theorem is given. We want to 
construct local weighted coordinates $x_a$ such that $Q$ is defined by the equations from Lemma \ref{lem:qincoord}. 
The proof (modeled after the construction of `privileged coordinates'  in Choi-Ponge \cite{choi:priv1,choi:priv2}, 
which in turn was motivated by Bella\"iche \cite{bel:tan}) is rather long, and will be given in the final section of this paper. 

\subsection{Application: Singular Lie filtrations} Using Theorem \ref{th:intrinsic}, we obtain new examples of weightings arising from \emph{singular Lie filtrations}. We will give a detailed discussion elsewhere \cite{loi:lie}, but here is a quick preview. A singular Lie filtration on a manifold $M$ is a filtration 
$\mf{X}_M=\H_{-r}\supset \cdots \supset \H_0$
of the sheaf of vector fields  by locally finitely generated $C^\infty_M$-submodules, compatible with brackets in the sense that $[\H_{-i},\H_{-j}]\subset \H_{-i-j}$. Examples of singular Lie filtrations abound:  they include  the much-studied \emph{regular} 
Lie filtrations, with $\H_{-i}$ the sheaves of sections of vector bundles (see references given in the introduction), but also many examples from sub-Riemannian geometry \cite{bel:tan}
or examples related to singular hypo-elliptic operators \cite{and:inprep}.  
 Let 
\[ \E\subset \mf{X}_{T_rM}\] 
be the $C^\infty_{T_rM}$-submodule spanned by  lifts $X^{(-i)}$ of vector fields $X\in \H_{-i}(U)$. The bracket condition implies that $\E$ is involutive, hence defines a singular foliation of $T_rM$ (in the sense of Androulidakis-Skandalis \cite{and:hol}). Suppose $Q\subset T_rM$ is a closed leaf. Then $Q$ satisfies the conditions of Theorem \ref{th:intrinsic}: Condition (iii) holds by construction, (ii) holds because the $\H_{-i}$ define a filtration, and (i) holds since $\H_{-r}=\mf{X}_M$. Hence, $Q$ defines a weighting of $M$ along $N=Q\cap M$ (which is itself a leaf of the singular foliation $\H_0$).

\subsection{Application: Multiplicative weightings} 
Viewing weightings in terms of subbundles of the $r$-th tangent bundle suggests notions of \emph{multiplicative weightings} on Lie groupoids and Lie algebroids. 
Given a Lie groupoid $G\rra M$ and a subgroupoid $H\rra M$, we say that a weighting of $G$ along $H$ is \emph{multiplicative} if the corresponding subbundle $Q\subset T_rG$ is a Lie subgroupoid of $T_rG\rra T_rM$. The units of this subgroupoid 
are then a graded subbundle $Q_N\subset T_rM$ defining a weighting of $M$ along $N$. We obtain groupoids
\[ \nu_\W(G,H)\rra \nu_\W(M,N),\ \ \delta_\W(G,H)\rra \delta_\W(M,N),\ \
\wh{\on{Bl}}^0_\W(G,H)\rra  \wh{\on{Bl}}_\W(M,N)\]
where  the weighted blow-up groupoid is defined as  the quotient of the \emph{restriction} of $\delta_\W(G,H)$ to ${\delta_\W(M,N)-(N\times\R)}$ by the action of $\R^\times$. (Recall that the restriction of a groupoid  to a subset of its units  is the set of arrows for which both source and target are in that subset.) See \cite{deb:blo,gua:sym} for blow-ups of Lie groupoids in the unweighted setting. 

Similarly, given a Lie algebroid $A\Rightarrow M$ with a Lie subalgebroid $B\Rightarrow N$, we can define \emph{infinitesimally multiplicative}  weightings of $A$ along $B$ by the condition that the corresponding graded subbundle $Q\subset T_rA$ is a Lie subalgebroid 
$Q\Rightarrow Q_N$
of  $T_rA\Rightarrow T_rM$. Given such a weighting, we obtain new Lie algebroids 
\[ \nu_\W(A,B)\Rightarrow \nu_\W(M,N),\ \ \delta_\W(A,B)\Rightarrow \delta_\W(M,N),\ \ \wh{\on{Bl}}^0_\W(A,B)\Rightarrow  \wh{\on{Bl}}_\W(M,N).\]

As an example, any order $r$ weighting of $M$ along $N$, defined by a graded subbundle $Q_N\subset T_rM$, 
the pair groupoid $\on{Pair}(M)=M\times M\rra M$ has a multiplicative weighting along $\on{Pair}(N)\rra N$, defined by the graded subbundle 
\[ \on{Pair}(Q_N)\subset \on{Pair}(T_rM)=T_r(\on{Pair}(M)).\]
Likewise, the tangent bundle $TM\Rightarrow M$ acquires an infinitesimally multiplicative weighting along $TN\Rightarrow N$, with the corresponding graded subbundle $TQ\subset T(T_rM)\cong T_r(TM)$.   One has an isomorphism $\nu_\W(TM,TN)\cong T\nu_\W(M,N)$, by an argument extending that in \cite[Appendix]{bur:spl}. 

\section{Further generalizations}
\subsection{Multi-weightings}\label{subsec:multiweightings}
The theory of weightings on manifolds extends to  \emph{multi-weightings},  with weight sequences $\mathbf{w}_1,\ldots,\mathbf{w}_n\in \Z_{\ge 0}^d$
for some $d>1$. We declare that a monomial $x^s=x_1^{s_1}\cdots x_n^{s_n}$ with $s=(s_1,\ldots,s_n)$ has multi-weight $s\cdot \mathbf{w}=\sum_a s_a \mathbf{w}_a$. For $U\subset \R^n$ open, we take $C^\infty(U)_{(\mathbf{i})}$ to be the ideal spanned by all monomials $x^s$ such that $s\cdot \mathbf{w}\ge \mathbf{i}$ (using the partial ordering given by the componentwise inequality). Generalizing the case $d=1$, this may be used as a local model to define multi-weightings on manifolds in terms of a filtration of the function sheaf by ideals $C^\infty_{M,(\mathbf{i})}$.  

\begin{example}
Let $N_1,\ldots,N_d\subset M$ be closed submanifolds, 	intersecting cleanly in the sense that 
every point of  $M$ is contained in a coordinate chart in which the submanifolds look like coordinate subspaces, that is, subspaces of $\R^n$ given by the vanishing of some of the coordinate functions.
Then there is a \emph{trivial multi-weighting}, where $C^\infty_{M,(i_1,\ldots,i_d)}$ consists of all smooth functions vanishing to order $i_\alpha$ on $N_\alpha$, for $\alpha=1,\ldots,d$. 
\end{example}

Let $e_\alpha\in \Z^d$ denote the $\alpha$-th standard basis vector. 
For a $d$-fold multi-weighting, we obtain (single) weightings of $M$ given 
$C^\infty_{M,(i e_\alpha)}$. The closed submanifolds 
\[ N_1,\ldots,N_d\subset M\]
defined by these weightings intersect cleanly, as defined above. (Using multi-weighted coordinates, these submanifolds look like coordinate subspaces.) 
Conversely, the collection of these weightings determines the multi-weighting via $C^\infty_{M,(\mathbf{i})}=\bigcap_\alpha C^\infty_{M,(i_\alpha e_\alpha)}$.  
We may thus regard a $d$-fold multi-weighting as a collection of $d$ weightings with some compatibility condition. 

One defines $d$-fold \emph{multi-graded bundles}  \cite{gra:hig,gra:gra}   as manifolds $E$ with smooth actions of the monoid $(\R^d,\cdot)=(\R,\cdot)\times \cdots \times (\R,\cdot)$. Put differently, $E$ comes equipped with $d$ commuting scalar multiplications $\kappa_{t_1}^1,\ldots,\kappa^d_{t_d}$. 
	The linear approximation $E_\lin=\nu(E,N)$ of a multi-graded bundle is a multi-graded \emph{vector bundle}. 
Multi-graded bundles have a canonical multi-weighting, described 
by the obvious generalization of \eqref{eq:weightingsgradedbundles}. 
An example of a multigraded bundle  is 
\[ T_{r_1,\ldots,r_d}M=J_0^{r_1,\ldots,r_d}(\R^d,M)=
\Hom(C^\infty(M),\AA_{r_1}\otimes \cdots \otimes \AA_{r_d}).\]
For instance, $T_{1,1}M\cong TTM$ is the double tangent bundle.

The construction of weighted normal bundle admits a straightforward generalization to multi-weightings. Given a multi-weighting of $M$ defined by a multi-filtration of $C^\infty(M)$, 
let $\on{gr}(C^\infty(M))$ be the associated multi-graded algebra, with summands
\[ C^\infty(M)_{(\mathbf{i})}/\sum_{\mathbf{j}>\mathbf{i}} C^\infty(M)_{(\mathbf{j})}.\]
The character spectrum 
\[\nu_\W(M,N_1,\ldots,N_d)=\Hom_\alg(\on{gr}(C^\infty(M)),\R)\]
is then a multi-graded bundle over the intersection $N=N_1\cap\cdots \cap N_d$. 
(Multi-graded coordinates on this bundle 
are constructed from multi-weighted coordinates on $M$, by a 
straightforward generalization from the $d=1$ case.) Likewise, we obtain $d$-fold deformation spaces
\[ \pi\colon \delta_\W(M,N_1,\ldots,N_d)\to \R^d.\] 

It admits an algebraic definition as the character spectrum of the $\Z^d$-graded Rees algebra of $C^\infty(M)$, consisting of  
of $d$-variable Laurent polynomials 
\[ \sum_{i_1,\ldots,i_d\in \Z} z_1^{-i_1}\cdots z_d^{-i_d}
f_{i_1,\ldots,i_d},\ \ \ f_{i_1,\ldots,i_d}\in C^\infty(M)_{(i_1,\ldots,i_d)}; \] 
the natural inclusion $\R[z_1,\ldots,z_d]\to \on{Rees}(C^\infty(M))$ defines the projection $\pi$. 
The fibers of the projection $\pi$ 
are themselves multi-weighted normal bundles (for smaller $d$). For instance, in the case of a bi-weighting,  
\[ \pi^{-1}(t_1,t_2)=\begin{cases}\nu_\W(M,N_1,N_2) & t_1=t_2=0\\
\nu_\W(M,N_1)& t_1=0,\ t_2\neq 0\\
\nu_\W(M,N_2)& t_1\neq 0,\ t_2=0\\
\nu(M,\emptyset)=M & t_1,t_2\neq 0\end{cases}\]

\begin{example}
	For a trivial multi-weighting along cleanly intersecting submanifolds $N_1,\ldots,N_d$, we obtain a $d$-fold vector bundle (called the $d$-fold normal bundle) 
	\[ 	\nu(M,N_1,\ldots,N_d)\to N,\]
	and a corresponding $d$-fold deformation space. See \cite{me:wei,pik:the} for $d=2$. 
	In \cite{pik:the}, it is shown that the double deformation space may also be obtained by iteration; the argument generalizes to $d>2$. 
\end{example}


Given a $d$-fold multi-weighting of $M$, let $r_1,\ldots,r_d$ be the orders of the resulting weightings along  the submanifolds $N_1,\ldots,N_d$, and $Q_1\subset T_{r_1}M,\ldots,Q_d\subset T_{r_d}M$ the corresponding graded subbundles.  By iterated tangent prolongation, these give rise to multi-graded subbundles $Q_\alpha'\subset T_{r_1,\ldots,r_d}M$ with clean intersection. 
The weighted normal bundle $\nu_\W(M,N_1,\ldots,N_d)$ is realized as a quotient of a multi-graded subbundle $Q=Q_1'\cap\cdots \cap Q_d'\subset T_{r_1,\ldots,r_d}M$, similar to the $d=1$ case. 


\begin{remark}
	To any multi-weighting  we may associate its \emph{total weighting} along the intersection $N=N_1\cap\cdots \cap N_d$, by
	taking $C^\infty(M)_{(j)}$ to be the sum of all $C^\infty(M)_{(\mathbf{i})}$ such that $|\mathbf{i}|\ge j$. 
	Note that a trivial multi-weighting may result in a non-trivial total weighting. Examples \ref{ex:ex}\eqref{it:nested} and \ref{ex:ex}\eqref{it:clean} are special cases of this construction (note that a nested sequence of submanifolds has clean 
	intersection). The weighted normal bundle for the total weighting is obtained by passing to the diagonal $(\R,\cdot)$-action as in Example \ref{ex:gradedbundles} \eqref{it:f} above. 
\end{remark}

\subsection{The holomorphic setting}\label{subsec:holo}
In this paper, we developed the theory of weightings  in the category of $C^\infty$-manifolds. Working with sheaves, the definition generalizes to complex manifolds, where it 
leads to a concept of holomorphic weighted normal bundle $\nu_\W(M,N)\to N$ with an associated holomorphic weighted 
deformation space $\delta_\W(M,N)$. While the main ideas of this extension are reasonably straightforward, there are some  aspects deserving special attention.  Note first that \cite{joz:monoid} describes the theory of graded bundles in the complex setting, using holomorphic actions of the monoid $(\C,\cdot)$ of complex numbers. 
One important difference to the $C^\infty$-context is Theorem \ref{th:gr} -- that is,  a holomorphic graded bundle is not, in general, globally isomorphic to its linear approximation.  

An example of a holomorphic graded bundle is the holomorphic $r$-th tangent bundle $T_rM$. This may be defined using holomorphic jets (see \cite{joz:monoid}), or algebraically using the holomorphic counterpart to definition \ref{def:algebraic}. A weighting is then characterized as a certain holomorphic graded subbundle $Q\subset T_rM$, and the weighted normal bundle is a quotient of the latter.  Holomorphic weighted Euler-like vector fields and the associated weighted tubular neighborhood embeddings only exist locally, in general -- in fact, the deformation space serves as a substitute for tubular neighborhood embeddings in some situations. 
Finally, in the complex setting there is no good analogue of the `spherical' weighted blow-up discussed in Section \ref{sec:weightedblowup}; on the other hand, the projective weighted blow-ups always acquire orbifold singularities.



\section{Proof of Theorem \ref{th:intrinsic}}\label{app:technical}

This  section is devoted to a proof of Theorem \ref{th:intrinsic}, characterizing the subbundles $Q\subset T_rM$ corresponding to weightings.  The approach is inspired by the construction of privileged coordinates in \cite{bel:tan,choi:priv1,choi:priv2}, and some of the techniques in \cite{hig:eul}.

\subsection{Filtration of vector fields}
Let $N\subset M$ be a closed submanifold. We denote by $\I\subset C^\infty_M$ the ideal sheaf of functions vanishing on $N$. 
Let $Q\subset T_rM$ be a graded subbundle along $N\subset M$. 	Define a filtration $C^\infty_M=C^\infty_{M,(0)}\supset 
C^\infty_{M,(1)}\supset\cdots
$,  by the prescription from Theorem \ref{th:q}: For $1\le i\le r$, \begin{equation}\label{eq:candidatefiltration} C^\infty(U)_{(i)}=\{f\in C^\infty(U)|\ \forall j< i\colon f^{(j)}	\mbox{ vanishes on}\ Q\cap T_rU\};\end{equation}
for $i>r$, 	we define $C^\infty(U)_{(i)}$ as sums of products of functions such that the filtration degrees add to $i$.
As usual, the filtration of the algebra of functions determines a filtration of vector fields 
$\mf{X}_M=\mf{X}_{M,(-r)}\supset \mf{X}_{M,(-r+1)}\supset \cdots$, where $X\in\mf{X}(U)$ has filtration degree $j$
if it takes $C^\infty(U)_{(i)}$ to $C^\infty(U)_{(i+j)}$. Note that $C^\infty(U)_{(1)}$ are the functions vanishing on $U\cap N$. 

Let us now impose the first two conditions from Theorem \ref{th:intrinsic}: thus $Q$ is invariant under the bundle action of $TM$
and its tangent bundle is $\AA_r$-invariant. Define a sheaves of vector fields by 
\begin{equation}\label{eq:ku} \K_{-i}(U)=\{X\in\mf{X}(U)|\ X^{(-i)}\mbox{ is tangent to } Q\cap T_rU\}.\end{equation} 
\begin{proposition}\label{lem:degofk}
The sheaves $\K_{-i}$ are $C^\infty_M$-submodules, with $[\K_{-i},\K_{-j}]\subset \K_{-i-j}$. They define a 
filtration of the sheaf of vector fields
\begin{equation}\label{eq:hfiltration} \mf{X}_M=\K_{-r}\supset \K_{-r+1}\supset \cdots \supset \K_0.\end{equation}
The submodule $\K_{-i}$ is contained in $\mf{X}_{M,(-i)}$, for $i=0,\ldots,r$. 	
\end{proposition}
\begin{proof}
 From $(fX)^{(-i)}=\sum_{j=0}^{r-i} f^{(j)} X^{(-i-j)}$ we see that the $\K_{-i}$ are $C^\infty_M$-submodules, and 
$[X^{(-i)},Y^{(-j)}]=[X,Y]^{(-i-j)}$ implies compatibility with Lie brackets. 
Furthermore, $\mf{X}_M=\K_{-r}$ by $TM$-invariance of $Q$, and $\K_{-i}\subset \K_{-i-1}$ by $\AA_r$-invariance of $TQ$.
It remains to show that $\K_{-i}\subset \mf{X}_{M,(-i)}$. Given $X\in \K_{-i}(U)$ and 
		$f\in C^\infty(U)_{(j)}$, we have to show  $Xf\in C^\infty(U)_{(j-i)}$.
		We may assume $j\ge i$.  Let $i\le \ell<j$. By \eqref{eq:liftdef}, 
		\[  (Xf)^{(\ell-i)}= X^{(-i)}f^{(\ell)}.\]
		But $f^{(\ell)}$ vanishes on $Q$ for $\ell<j$, and 
	 $X^{(-i)}$ is tangent to 
		$Q\cap T_rU$. It hence follows that $(Xf)^{(\ell-i)}$ vanishes on $Q\cap T_rY$. This means that $Xf$ has filtration degree $j-i$. 	
	\end{proof}
	\begin{remark}
		Proposition \ref{prop:gq} shows that for graded subbundles $Q\subset T_rM$ defined by a weighting, the submodules
		$\K_{-i}$ coincide with $\mf{X}_{M,(-i)}$ for $i=0,\ldots,r$. On the other hand, 
		 Example \ref{ex:counterexample} shows that in general, the inclusion  $\K_{-i}\subseteq\mf{X}_{M,(-i)}$ is strict.  
	\end{remark}
\begin{lemma}\label{lem:takes}
	For $0\le j\le i\le r$,  multiplication by $C^\infty_{M,(j)}$ takes $\K_{-i}$ to $\K_{-i+j}$. 
\end{lemma}
\begin{proof}
	Let $g\in C^\infty(U)_{(j)}$ and $X\in \K_{-i}(U)$. We have to show that $(gX)^{(-i+j)}=\sum_\ell  g^{(\ell)}X^{(-i-\ell+j)}$
	is tangent to $Q\cap T_rU$. For $\ell<j$ the function $g^{(\ell)}$ vanishes 
on $Q\cap T_rU$, while for $\ell\ge j$ we have that $-i-\ell+j\le -i$, which implies that $ X^{(-i-\ell+j)}$ is tangent to $Q\cap T_rU$. 
\end{proof}

	Let $\on{DO}_M$ be the sheaf of scalar differential operators on $M$, and $\on{DO}^p_M$ the subsheaf of differential operators of order at most $p$.  (That is, $\on{DO}^p(U)$ are scalar differential operators of order $\le p$ acting on $C^\infty(U)$.) The filtration of vector fields, given by the submodules $\K_{-i}$, extends to a filtration on each $\on{DO}^p(M)$. 
	
	\begin{definition}\label{def:qweight}
	We say that $D\in \on{DO}(U)$ \emph{has $Q$-weight $\ell\le 0$} if it can be written as a linear combination of operators 
	\[ X_p\cdots X_1,\]
	with 
	$X_{\nu}\in \K_{-j_\nu}(U)$ and $-(j_1+\ldots+j_p)=\ell$. Let $\on{DO}^{p}(U)_\ell$ be the space of differential operators of order $\le p$ and $Q$-weight $\ell$. 
	\end{definition}
   	For fixed $p$, the filtration starts by $Q$-weight starts in degree $-rp$: 
   	\[ \on{DO}^{p}_{M}=\on{DO}^{p}_{M,-rp}\supset \cdots \supset \on{DO}^{p}_{M,0}\supset 0.\]
   	As a consequence of Lemma \ref{lem:degofk}, and a straightforward induction, the filtration is compatible with products: If $D_1,D_2$ have $Q$-weights $\ell_1,\ell_2$, respectively, then $D_1\circ D_2$ has $Q$-weight $\ell_1+\ell_2$. 	Lemma \ref{lem:degofk} also shows that differential operators of $Q$-weight $\ell$ act on functions as operators of 
	filtration degree $\ell$. 
	
	\subsection{Consequences of transitivity}
In terms of the submodules $\K_{-i}$, property (iii) from Theorem \ref{th:intrinsic} says that the tangent bundle of $Q\cap T_rU$
is spanned by all $X^{(-i)}$ with $i=0,\ldots,r$ and $X\in \K_{-i}(U)$. In other words, the action of the Lie subalgebra 
\[ \big\{\sum_{i=0}^r X_i\otimes \epsilon^i\colon X_i\in\K_{-i}(U)\big\}\subset \mf{X}(U)\otimes \AA_r\]
given by $\sum_{i=0}^r X_i\otimes \epsilon^i\mapsto \sum_{i=0}^r X_i^{(-i)}$ is infinitesimally transitive. 
Let us assume this transitivity  property from now on. 	
	Define $\wt{F}_{-i}=Q_\lin^{-i}$ for $i=1,\ldots,r$ and $\wt{F}_0=TN$. Thus
	\[ TQ|_N=\wt{F}_0\oplus \cdots \oplus \wt{F}_{-r},\]
	as a subbundle of $T(T_rM)=TM^{\oplus(r+1)}$. Transitivity of the action implies that $\wt{F}_{-i}$ is spanned by the restriction of 
	$\K_{-i}$ to $N$. We hence obtain a filtration by subbundles, 
	\begin{equation}\label{eq:filtrationofTMN}
	TM|_N= \wt{F}_{-r}\supset \wt{F}_{-r+1}\supset \cdots \wt{F}_{-1}\supset \wt{F}_0=TN.\end{equation}
	Letting $k_i=\on{rank}(\wt{F}_{-i})$ we obtain  a unique  weight sequence $0\le w_1\le \cdots\le w_n\le r$ 
	such that $w_a<i \Leftrightarrow a< k_i$. 
	We may characterize the filtration of functions in terms of the filtration of vector fields by $\K_{-i}\subset \mf{X}_M$.

	\begin{lemma}\label{lem:crit1} Let $1\le i\le r$. 	A function $f\in C^\infty(U)$ has filtration degree $i$ if and only if 
		\[ Df|_N=0\]
		for all $D\in \on{DO}(U)$ of $Q$-weight $\ell$, with $\ell+i>0$. 
		
	\end{lemma}
	\begin{proof}		
		Suppose that $f\in C^\infty(U)$ has filtration degree $i$.  
		To show $Df|_N=0$ for all 
		differential operators of $Q$-weight $\ell>-i$,  it 
		suffices to check for  $D=X_{p}\cdots X_{1}$ with 
	 $X_{\nu}\in \K_{-j_\nu}$ where $-(j_1+\ldots+j_p)>-i$. 
		By Lemma \ref{lem:degofk}, vector fields in $\K_{-j}$ lower the filtration degree on functions by $j$. 
		Hence, if $f$ has filtration degree $i$, then $D f$
		has filtration degree $i-(j_1+\ldots+j_p)>0$, and so $Df|_N=0$. 
		
		For the converse, suppose that $f\in C^\infty(U)$ satisfies the conditions of the lemma, for some given $i$.  To show that $f$ has filtration degree $i$, we need to show that $f^{(j)}$ vanishes on $Q$ for all $j<i$. 
	 By the transitivity assumption, the action of the nilpotent Lie subalgebra  
	 \[ \big\{\sum_{\ell=1}^r X_\ell\otimes \epsilon^\ell\colon X_i\in \K_{-i}(U)\big\}\] by vertical vector fields $\sum_{\ell=1}^r X_\ell^{(-\ell)}$ is 
	 infinitesimally transitive on the fibers of $Q\cap T^rU\to N\cap U$. 		The elements $\exp(\sum_{\ell=1}^r X_\ell\otimes \epsilon^l),\ X_\ell\in \K_{-\ell}(U)$
		of the corresponding nilpotent group act by the time-one flow of $\sum_{\ell=1}^r X_\ell^{(-\ell)}$. 
		By the Baker-Campbell-Hausdorff formula, we may write
		\[ \exp(\sum_{\ell=1}^r X_\ell\otimes \epsilon^l)=\exp(\sum_{\ell=j+1}^r X_\ell' \otimes\epsilon^\ell)\exp(\sum_{\ell=1}^j X_\ell \otimes\epsilon^\ell)\]
		with new elements $X_\ell'\in \K_{-\ell}(U)$. The action of $\exp(\sum_{\ell=j+1}^r X_\ell' \otimes\epsilon^\ell)$ on $T_rU$ preserves the fibers of $p^r_j\colon T_rU\to T_jU$, and $f^{(j)}$ is constant along those fibers. Hence, to show that 
		$f^{(j)}$ vanishes on $Q\cap T_rU$, it suffices to show that $(\Phi_{t}^Z)_*f^{(j)}|_{N\cap U}=0$ for all vector fields of the form
		\begin{equation}\label{eq:z}
		Z=\sum_{\ell=1}^j X_\ell^{(-\ell)},\ \ \ X_\ell\in \K_{-\ell}(U),\end{equation}		
		with $\Phi_t^Z$ denoting the flow of $Z\in \mf{X}(T_rU)$. 
		Note that the flow $\Phi_t^Z$ preserves the 
		spaces of polynomials of  any given degree $\le d$, since the Lie derivative $\L_{Z}$ does. Hence 
		$(\Phi_{t}^Z)_*$ is given on fiberwise polynomial functions, such as $f^{(j)}$, by  a truncated  exponential series
		\[ (\Phi_{t}^Z)_*f^{(j)}=\exp(t \L_{Z})f^{(j)}=\sum_{p=0}^j \f{t^p}{p!} (\L_{Z})^p f^{(j)}.\]
		We have 
		\begin{eqnarray*}
			\L_{Z}f^{(j)}&=&\sum_{j_1=1}^j  (\L_{X_{j_1}}f)^{(j-j_1)},\\
			(\L_{Z})^2f^{(j)}&=&\sum_{j_1=1}^j\sum_{j_2=1}^{j-j_1} 
			(\L_{X_{j_2}}\L_{X_{j_1}}f)^{(j-j_1-j_2)},
		\end{eqnarray*}	
		and so on. 	Upon restriction to $U\subset T_rU$, only the terms with $j_1+j_2+\ldots+j_p=j$ remain: 
		\[(\L_{Z})^p f^{(j)}|_U=\sum_{j_1+\ldots+j_p=j} \L_{X_{j_p}}\cdots \L_{X_{j_1}} f.\]
		By assumption, this expression vanishes along $N\cap U$, and so $f^{(j)}$ vanishes on $Q\cap T_rU$. 
	\end{proof}
	
	\subsection{Local frames}
	By an \emph{adapted local frame} over an open subset $U\subset M$, we mean a  frame $V_1,\ldots,V_n\in\mf{X}(U)$ of $TM|_U$, such that for all $i\le r$, 
	\[ 
	V_1,\ldots,V_{k_i}\in \K_{-i}(U),
	\]
	and such that the restrictions $V_a|_N$ for $a\le k_0$ commute (recall that vector fields in $\K_0(U)$ are tangent to $N\cap U$).
	Given $m\in M$, one may construct an adapted local frame over a neighborhood $U$ of $m$: Start with a frame of $\ti{F}_0=TN$, given by a collection of commuting vector fields on $N$. Extend to a 
	local frame for $TM|_N\to N$ near $m$, adapted to the filtration by subbundles $\ti{F}_{-i}$. 
	Then  use $\ti{F}_{-i}=\K_{-i}|_N$  to extend these sections to a local frame for $TM\to M$ near $m$. 
	
	The lifts 	$V_a^{(-i)}\in \mf{X}(T_rU)$ for $i=0,\ldots,r$ and $a=1,\ldots,n$ are a frame for $T(T_rU)$. Those satisfying the extra condition $a\le k_i$ (i.e., $w_a\le i$)  restrict to a frame for $T(Q\cap T_rU)$: they are tangent to $Q$ by definition 
	of $\K_{-i}$, and they span the tangent bundle by dimension count. Every $D\in \on{DO}(U)$ can be uniquely written
	as a finite sum 
	\begin{equation}\label{eq:standardform} D=\sum_{|s|\le p} f_s V^s,\end{equation}
	using multi-index notation  $s=(s_1,\ldots,s_n),\ |s|=\sum_a s_a$ with 
	$V^s=V_1^{s_1}\cdots V_n^{s_n}$. 
	
	\begin{lemma}
		If $D\in \on{DO}(U)$ has $Q$-weight $\ell$, then the functions $f_s$ defined by the standard form 
		\eqref{eq:standardform}  have filtration degree $\ell+w\cdot s$.
	\end{lemma}
	\begin{proof}
		Let $D$ be a differential operator of $Q$-weight $\ell\le 0$. If $D$ has order $0$ (so that it is given by a function, acting by multiplication) the statement is obvious: Every function has filtration degree $0$, hence also filtration degree $\ell\le 0$.
		Consider next the case that $D$ is a vector field $X\in \K_{-i}(U)$ with $\ell=-i$.  
		Writing $X=\sum_{a=1}^n f_a V_a$ we want to show that the coefficient functions $f_a\in C^\infty(U)$ have filtration degree $w_a-i$. 
		By definition of $\K_{-i}$, the vector field 
		\[ X^{(-i)}=\sum_{a=1}^{n} \sum_{j=0}^{r-i}  f_a^{(j)} V_a^{(-i-j)}\]
		is tangent to $Q\cap T_rU$. But $V_a^{(-i-j)}$ 
		belongs to the frame for $T(Q\cap T_rU)$ if and only if $w_a\le i+j$. 
		Hence, $X^{(-i)}$ is tangent to $Q\cap T_rU$ if and only if 
		$f_a^{(j)}|_Q=0$ for $w_a>i+j$, which is the case if and only if $f_a$ has filtration degree $w_a-i$. This proves the claim for $D$ of order $\le 1$.
		
		To prove the general case, it hence suffices to show that if the statement holds for differential operators $D_1,D_2$ of $Q$-weights $\ell_1,\ell_2$, then it also holds for the product $D=D_1\circ D_2$, as a differential operator of 
		$Q$-weight $\ell=\ell_1+\ell_2$. 
		This involves re-arranging the terms in 
		$D_1 \circ D_2$ to bring them to the standard form $\sum_{|s|\le p} f_s V^s$. It is enough to check this on generators for the algebra $\on{DO}(U)$. If $a>b$ (so that $V_a\circ V_b$ is in the `wrong' order),  we have 
		\[ V_a\circ V_b=V_b\circ V_a+[V_a,V_b].\]
		Here $X=[V_a,V_b]\in \K_{-w_a-w_b}$, and (as shown above) the coefficients $f_{ab}^c$ in the bracket relation  $[V_a,V_b]=\sum_c f_{ab}^c V_c$ have filtration degree $w_c-w_a-w_b$. On the other hand, if 
		$f\in C^\infty(U)$ and $V_a\in \K_{-w_a}(U)$, consider $V_a\circ f$ as a product of differential operators of $Q$-weights $\ell_1=-w_a$ 
		and $\ell_2=0$.  Its standard form reads as 
		\[ V_a\circ f=f  V_a+V_a(f).\]
		As noted above, for functions (such as $V_a(f)$) the claim is obvious. On the other hand, 
		since $\K_{-w_a}(U)$ is a $C^\infty(U)$-module,  the first term is a differential operator of $Q$-weight 
		$-w_a$, and it it has the required form since the coefficient function $f$ has filtration degree $\ell+w\cdot s=
		-w_a+w_a =0$.   
	\end{proof}

	This allows us to reformulate Lemma  \ref{lem:crit1} as follows: 
	
	\begin{lemma}\label{lem:vslemma}
		A function $f\in C^\infty(U)$ has filtration degree $i$ if and only if 
		\[ (V^s f)|_N=0\]
		for all multi-indices $s$ with $w\cdot s<i$. (It suffices to verify this condition for multi-indices with $s_a=0$ for $a\le k_0$.) 
	\end{lemma}
	\begin{proof}
		The condition is necessary, since the differential operator $V^s$ has $Q$-weight $-w\cdot s$, and hence has filtration degree $-w\cdot s$ as an operator on functions. Conversely, suppose that the condition is satisfied.  
		Lemma \ref{lem:crit1} shows that $f$ has filtration degree $i$ if and only if $Df|_N=0$ for all differential operators 
		$D$ of $Q$-weight $\ell>-i$. Write any such $D$ in its standard form \eqref{eq:standardform} where $f_s$ has filtration degree $\ell+w\cdot s$. For coefficients $s$ with 
		$\ell+w\cdot s>0$, already $f_s$ vanishes on $N$. For coefficients with $\ell+w\cdot s\le 0$, the assumption guarantees that $(V^sf)|_N=0$. Hence $(Df)|_N=0$ as required. 
	\end{proof}
	
	\subsection{Adapted coordinates}
	After these preparations, we get to the last stage of the proof, using the frame $\{V_a\}$ 
	on a neighborhood $U$ of a given point $m\in N$ to construct adapted coordinates, so that 
	$Q\cap T_rU$ is given by the equations from Lemma \ref{lem:qincoord}. We observe that this is equivalent to $x_a$ having the expected filtration degrees (using the filtration \eqref{eq:candidatefiltration} on functions):

	\begin{lemma}
		Suppose $U\subset M$ is an  open neighborhood of $m$, with coordinate functions
		$x_a\in C^\infty(U)$, such that $x_a$ has filtration degree $w_a$. Then the manifold $Q\cap T_rU$ is given by the equations 
		\[ x_a^{(i)}=0\ \ \mbox{ for }\ \ i=0,\ldots,r,\ a>k_i.\] 
		In particular, $Q\cap T_rU$ corresponds to a weighting of $(U,N\cap U)$.
	\end{lemma}
	\begin{proof}
		By definition, the coordinate function $x_a$ has filtration degree $w_a$ if and only if 
		the function $x_a^{(i)}$  vanishes on $Q\cap T_rU$ for all $a>k_i$, i.e.  $w_a>i$.  Since 
		$\dim Q=k_0+\ldots+k_r$, it follows that $Q\cap T_rU$ is given exactly by the vanishing of these functions. 
	\end{proof}

	\begin{lemma}\label{lem:finally}
		After replacing $U$ with a smaller neighborhood of $m$, if necessary, there are local coordinates $x_1,\ldots,x_n$ on $U$ such that each $x_a$ has filtration degree $w_a$, and 
		\[ (V_a\,x_b)|_{N\cap U}=\delta_{ab}.\] 
	\end{lemma}
	\begin{proof}
		Since 	the restrictions  $V_a|_N,\ 1\le a\le k_0$ are commuting vector fields on $N$ (near $m$), we can choose functions $y_a\in C^\infty(U)$ for $1\le a\le k_0$ such that 
		\begin{equation}\label{eq:deltaab}(V_a\, y_b)|_{N\cap U}=\delta^a_b
		\end{equation} 
		for $a,b\le k_0$. Complete to a collection of functions $y_a\in C^\infty(U),\ 1\le a\le n$, with $y_a|_{N\cap U}=0$ for $a>k_0$, 
		such that \eqref{eq:deltaab} holds for all 
		$1\le a,b\le n$. Our goal is to replace the $y_a$ with 
		new coordinates $x_a$, satisfying analogous properties, such that the $x_a$  have filtration degrees $w_a$. 
		For $a\le k_0$ we will put $x_a=y_a$. For $a>k_0$, we will
	    look for a coordinate change of the form  
		\[ x_a=y_a+\sum  \chi_{au}\,  y^u,\ \ \ \ \ a>k_0 \]
		(using multi-index notation  $y^u=y_1^{u_1}\cdots y_n^{u_n}$), 
		where the sum is over multi-indices with  $|u|=\sum_a u_a\ge 2$, with $u_b=0$ for $b\le k_0$, and 
		such that the coefficients $\chi_{au}\in C^\infty(U)$ depend only on the coordinates $y_1,\ldots,y_{k_0}$. 
		We will also impose the additional condition 
		$u\cdot w<w_a$ on the multi-indices, since terms with $u\cdot w\ge w_a$ would have filtration degree $\ge w_a$.
    Note that a coordinate change of this form adds terms quadratic and higher in the `normal coordinates' $y_a,\ a>k_0$; hence it will retain the property $(V_a x_b)|_N=\delta^a_b$. 
		
		By Lemma \ref{lem:vslemma}, the coordinate function $x_a$ has filtration degree $w_a$ if and only if $(V^s x_a)|_N=0$ for all multi-indices $s$ with $w\cdot s<w_a$. Again, we need only consider multi-indices where $s_a=0$ for $a\le k_0$. 
		For any such multi-index, this gives the condition 
		\[ 0=(V^s x_a)|_N=(V^s y_a)|_N+\sum_u \big(V^s(\chi_{au}\,  y^u)\big)\big|_N. \]
		If $|u|>|s|$ then  $V^s\big(\chi_{au}\,  y^u\big)$ is a polynomial of degree $>0$ in the normal coordinates, 
		and so vanishes on $N$. For $|u|=|s|$, we obtain a non-zero term only when $u=s$, and by applying all derivatives to $y^u$. Letting $c_s=V^s(y^s)|_N$ (a positive constant), we therefore obtain the condition
		\[ 0=(V^s y_a)|_N+c_s \chi_{as}|_N+\sum_{u\colon |u|<|s|} \big(V^s(\chi_{au}\,  y^u)\big)\big|_N.\]
		That is, 
		\[ \chi_{as}=-\f{1}{c_s}\Big(V^s y_a|_N+\sum_{u\colon |u|<|s|} \big(V^s(\chi_{au}\,  y^u)\big)\big|_N\Big).\]
		This gives a recursive formula for the desired coordinate change. 
	\end{proof}
	This completes the proof of Theorem \ref{th:intrinsic}, giving a characterization of weightings in terms of a graded subbundle $Q\subset T_rM$. 	

\bibliographystyle{amsplain} 

\def\cprime{$'$} \def\polhk#1{\setbox0=\hbox{#1}{\ooalign{\hidewidth
			\lower1.5ex\hbox{`}\hidewidth\crcr\unhbox0}}} \def\cprime{$'$}
\def\cprime{$'$} \def\cprime{$'$} \def\cprime{$'$} \def\cprime{$'$}
\def\polhk#1{\setbox0=\hbox{#1}{\ooalign{\hidewidth
			\lower1.5ex\hbox{`}\hidewidth\crcr\unhbox0}}} \def\cprime{$'$}
\def\cprime{$'$} \def\cprime{$'$} \def\cprime{$'$} \def\cprime{$'$}
\providecommand{\bysame}{\leavevmode\hbox to3em{\hrulefill}\thinspace}
\providecommand{\MR}{\relax\ifhmode\unskip\space\fi MR }
\providecommand{\MRhref}[2]{%
	\href{http://www.ams.org/mathscinet-getitem?mr=#1}{#2}
}
\providecommand{\href}[2]{#2}

\end{document}